\documentclass[12pt]{iopjournal}
\usepackage[T1]{fontenc}
\usepackage{lmodern}
\usepackage{amsmath,amssymb,amsthm,mathtools}
\usepackage[a4paper,margin=25mm,footskip=12mm]{geometry}
\usepackage{microtype}
\usepackage{ragged2e}
\usepackage{bookmark}
\usepackage{needspace}
\IfFileExists{xurl.sty}{\usepackage{xurl}}{}
\hypersetup{hidelinks,hypertexnames=false,
  pdftitle={Existence and multiplicity of radial solutions for higher-order logarithmic and infinite-mass Schrodinger equations},
  pdfauthor={Chen Huang, Rong Li, Junren Luo and Xian Zhang},
  pdfsubject={Schrodinger equation; higher-order logarithmic nonlinearity; infinite mass; radial potential; perturbation method; Pohozaev identity}}
\allowdisplaybreaks[2]
\numberwithin{equation}{section}
\newtheorem{theorem}{Theorem}[section]
\newtheorem{lemma}{Lemma}[section]

\theoremstyle{remark}
\newtheorem{remark}{Remark}[section]
\newcommand{\R}{\mathbb R}
\newcommand{\norm}[1]{\lVert #1\rVert}
\newcommand{\pair}[2]{\langle #1,#2\rangle}
\newcommand{\dx}{\,\mathrm{d}x}
\newcommand{\rad}{\mathrm{rad}}
\newcommand{\weakto}{\rightharpoonup}
\DeclareMathOperator{\spanof}{span}

\renewcommand{\articletype}[1]{\noindent{\small\sffamily #1}\par\medskip}
\renewcommand{\title}[1]{{\Large\bfseries\raggedright #1\par}\medskip}
\makeatletter
\renewcommand{\author}[1]{\noindent\if@anonymous
  Author names removed for review\else #1\fi\par\medskip}
\renewcommand{\affil}[1]{\if@anonymous\else
  {\small\noindent #1\par}\fi}
\renewcommand{\email}[1]{\if@anonymous\else
  \smallskip{\small\noindent\textbf{E-mail:} #1\par}\fi}
\makeatother
\renewcommand{\keywords}[1]{\smallskip
  {\small\noindent\textbf{Keywords:} #1\par}}
\renewenvironment{abstract}{\par\bigskip
  \noindent\textbf{Abstract}\par\smallskip
  \normalfont\normalsize\justifying\noindent\ignorespaces}{\par\medskip}

\begin{document}
\justifying
\articletype{}
\title{Existence and multiplicity of radial solutions for higher order logarithmic and infinite mass Schr\"odinger equations}
\author{Chen Huang$^1$, Rong Li$^2$, Junren Luo$^1$ and Xian Zhang$^{3,\ast}$}
\affil{$^1$School of Mathematics, Faculty of Basic Sciences, University of Shanghai for Science and Technology,
Shanghai 200093, China}
\affil{$^2$Department of Electronic and Information Engineering, School of Engineering, Westlake University, Hangzhou, Zhejiang 310030, China}
\affil{$^3$College of Economics and Management, Zhejiang University of Water Resources and Electric Power, Hangzhou, Zhejiang 310018, China}
\affil{$^\ast$Author to whom any correspondence should be addressed.}
\email{\href{mailto:chenhuangmath111@163.com}{chenhuangmath111@163.com} (Chen Huang); \href{mailto:lirong13@westlake.edu.cn}{lirong13@westlake.edu.cn} (Rong Li); \href{mailto:junrenluo@usst.edu.cn}{junrenluo@usst.edu.cn} (Junren Luo); \href{mailto:xianzhangmath@163.com}{xianzhangmath@163.com} (Xian Zhang)}

\begin{abstract}
We study radial solutions of the Schr\"odinger equation
$-\Delta u+V(x)u=g(u)$ in $\mathbb{R}^N$, $N\ge3$, where $V$ is radial
and $g(t)/t\to-\infty$ as $t\to0$. In this infinite-mass case the energy
functional is not finite on $H^1(\mathbb{R}^N)$. First let
$g(t)=t|\log|t||^{\alpha-1}\log|t|$ with $\alpha>1$. We only assume that
$V$ is continuous and bounded from below. So $V$ may change sign, and it
need not be coercive or have a limit at infinity. We obtain a nontrivial
radial solution and a sequence of radial solutions whose energies tend
to infinity. To do this, we add the term $\lambda|u|^{p-2}u$ with
$1<p<2$, find critical points of the perturbed functional by minimax
methods, and let $\lambda\to0$. The infinite-mass structure gives the
compactness in this limit. Next let $g$ be odd, with a power bound near
zero and at most critical growth at infinity. Let $V$ be bounded, with a
smallness condition on $(x\cdot\nabla V)_+$. We combine the
Pohozaev identity with an augmented symmetric minimax method and obtain
infinitely many pairs of radial solutions. A negative critical power
term is allowed.
\end{abstract}

\section{Introduction}

Logarithmic nonlinearities appear in several models of quantum mechanics
and nonlinear optics. Bia\l ynicki-Birula and Mycielski \cite{BBM}
introduced the equation
\begin{equation*}
 \mathrm{i}\hbar\partial_t\Psi
 =-\frac{\hbar^2}{2m}\Delta\Psi+V(x)\Psi
       -b\Psi\log\!\left(\frac{|\Psi|^2}{\rho_*}\right),
       \qquad \rho_*>0.
\end{equation*}
With the logarithmic term, noninteracting subsystems stay separable, and
the equation has Gaussian solitary waves. Logarithmic equations also
model self-sustained quantum Bose liquids \cite{AZ}. Logarithmically
saturable media appear in the study of incoherent optical solitons
\cite{Buljan}. In these models one looks for localized stationary states,
often without an external trapping potential. In all of them, however,
only the first power of the logarithm occurs.

Higher powers of the logarithm come from nonadditive entropies. The
Boltzmann--Gibbs entropy is additive for statistically independent
subsystems. This reflects the identity
$\log(\rho_1\rho_2)=\log\rho_1+\log\rho_2$. For systems with strong
correlations, long-range interactions or multifractal structures,
nonadditive entropies have been proposed; see
\cite{Tsallis1988,Tsallis2022}. Using fractional calculus, Ubriaco
\cite{Ubriaco} introduced entropies with densities $\rho[-\log\rho]^q$.
Stable equilibria minimize the free energy, and the dynamics is driven by
its variational derivative. So these entropies lead to nonlinearities
with higher powers of the logarithm. Al-zaleq and Alzaleq \cite{Alzaleq}
found exact solitary and periodic traveling waves for several logarithmic
nonlinearities. Li, Cui and Ye \cite{LCY} derived a higher-order
logarithmic Schr\"odinger equation from a nonadditive entropy. They
studied its traveling waves and spiral dynamics. In this paper we study
the existence and multiplicity of finite-energy stationary states.

We first consider
\begin{equation}\label{eq:original}
 -\Delta u+V(x)u=f(u)\quad\text{in }\R^N\quad(N\ge3),
\end{equation}
where
\begin{equation}\label{eq:fF}
 f(t)=\begin{cases}
 t|\log|t||^{\alpha-1}\log|t|,&t\ne0,\\
 0,&t=0,
 \end{cases}
 \qquad \alpha>1.
\end{equation}
For $\alpha=1$, \eqref{eq:original} is the logarithmic equation
$-\Delta u+V(x)u=u\log|u|$. The energy functional of \eqref{eq:original}
is
\begin{equation}\label{eq:energy0}
 J_0(u)=\frac12\int_{\R^N}\bigl(|\nabla u|^2+V(x)u^2\bigr)\dx
                  -\int_{\R^N}F(u)\dx,\qquad
 F(t)=\int_0^tf(s)\,\mathrm{d}s.
\end{equation}
Since $f(t)/t\to-\infty$ as $t\to0$, we have $-F(t)/t^2\to+\infty$ as
$t\to0$. As a result, $J_0(u)=+\infty$ for some $u\in H^1(\R^N)$. So
$J_0$ is not finite on $H^1(\R^N)$, and it is not of class $C^1$ there.

The case $\alpha=1$ has been studied extensively. Cazenave
\cite{Cazenave} worked in a Banach space with a Luxemburg-type norm to
overcome the lack of smoothness. In \cite{DMS}, d'Avenia, Montefusco and
Squassina used nonsmooth critical point theory. They obtained infinitely
many radial solutions for constant potentials. Squassina and Szulkin
\cite{SS} treated periodic potentials. They wrote the functional as the
sum of a $C^1$ part and a convex lower semicontinuous part. Ji and
Szulkin~\cite{JiSzulkin2016} obtained infinitely many solutions for
coercive potentials that may change sign. Zhang and
Zhang~\cite{ZhangZhang2022} proved existence, uniqueness and
nondegeneracy of positive solutions for radial potentials that may be
singular or tend to $-\infty$ at infinity. Wang and Zhang \cite{WZ}
showed that solutions of power-law equations converge to solutions of
the logarithmic equation. Huang and Yang \cite{HY} considered positive
potentials tending to infinity. They added $\lambda|u|^{p-2}u$ with
$1<p<2$ and did not truncate the logarithmic term. The same perturbation
was used in \cite{HYZ} to obtain infinitely many sign-changing
solutions.

For $\alpha>1$, An and Fang \cite{AF} studied ground states for an autonomous Schr\"{o}dinger equation with a power of the logarithm. In a recent preprint, Ou, Zhang and Li \cite{OZL} considered an equation
with both first- and higher-order logarithmic terms. They studied
nonnegative bound states for an asymptotically constant potential. In the
present paper we obtain infinitely many radial solutions of
\eqref{eq:original} for continuous radial potentials bounded from below.
We do not assume that $V$ is positive or monotone, or that it has a
limit at infinity.

\medskip
\noindent\textbf{Notation.} For $1\le s<\infty$,
$|u|_s=(\int_{\R^N}|u|^s\dx)^{1/s}$ is the norm of $L^s(\R^N)$, and
$|u|_\infty$ is the norm of $L^\infty(\R^N)$. We write
$\norm{u}^2=\int_{\R^N}(|\nabla u|^2+u^2)\dx$ for the norm of
$H^1(\R^N)$, $2^*=2N/(N-2)$, $a_\pm=\max\{\pm a,0\}$ and
$B_R=\{x\in\R^N:|x|<R\}$. We denote by $H^1_{\rad}(\R^N)$ and
$C_{c,\rad}^\infty(\R^N)$ the radial functions in $H^1(\R^N)$ and in
$C_c^\infty(\R^N)$. The letter $C$ denotes a positive constant that may
change from line to line.
\medskip

Our assumption on the potential is
\begin{equation}\tag{V1}\label{V1}
 V\in C(\R^N,\R),\quad V(x)=V(|x|)\quad\text{and}\quad
 V(x)\ge -M_0\quad\text{for }x\in\R^N,
\end{equation}
where $M_0\ge0$ is a constant. Let
\begin{equation*}
 \begin{split}
 E_V&=\left\{u\in H^1_{\rad}(\R^N):
                    \int_{\R^N}V_+(x)u^2\dx<\infty\right\},\\
 \norm{u}_V^2&=\int_{\R^N}
                  \bigl(|\nabla u|^2+(1+V_+(x))u^2\bigr)\dx.
 \end{split}
\end{equation*}
We say that $u\in E_V$ is a radial finite-energy weak solution of
\eqref{eq:original} if
\begin{equation}\label{eq:weak}
 \begin{gathered}
 \int_{\R^N}|u|^2|\log|u||^\alpha\dx<\infty,\\
 \int_{\R^N}\bigl(\nabla u\cdot\nabla\varphi+V(x)u\varphi\bigr)\dx
       =\int_{\R^N}f(u)\varphi\dx
       \quad\text{for all }\varphi\in C_c^\infty(\R^N).
 \end{gathered}
\end{equation}

We use a perturbation method. Fix $1<p<2$ and let
\begin{equation}\label{eq:X}
 X=E_V\cap L^p(\R^N),\qquad \norm{u}_X=\norm{u}_V+|u|_p.
\end{equation}
For $0\le\lambda\le1$, define
\begin{equation}\label{eq:energy}
 J_\lambda(u)=J_0(u)+\frac\lambda p\int_{\R^N}|u|^p\dx,
       \qquad u\in X.
\end{equation}
For $\lambda>0$, critical points of $J_\lambda$ solve
\eqref{eq:original} with the extra term $\lambda|u|^{p-2}u$ on the
left-hand side. The functional $J_\lambda$ is of class $C^1$ on $X$
(Section~\ref{sec:existence}). We find critical
points of $J_\lambda$ by minimax methods. We bound them uniformly in
$\lambda$. Then we let $\lambda\to0$ and obtain solutions of
\eqref{eq:original}. For quasilinear equations, this idea was used by
Liu, Liu and Wang \cite{LLW} with a $4$-Laplacian perturbation.

Our main results for \eqref{eq:original} are as follows.

\begin{theorem}\label{thm:existence}
Assume $N\ge3$, $\alpha>1$ and \eqref{V1}. Then \eqref{eq:original}
has a nontrivial radial finite-energy weak solution $u_0\in E_V$ with
\begin{equation*}
 J_0(u_0)=\lim_{\lambda\to0^+}c_\lambda>0,
\end{equation*}
where $c_\lambda$ is the mountain pass level of $J_\lambda$ defined in
\eqref{eq:MP-level}. Moreover, there are $\lambda_n\to0^+$ and mountain
pass critical points $u_n$ of $J_{\lambda_n}$ such that $u_n\to u_0$ in
$E_V$.
\end{theorem}

\begin{theorem}\label{thm:multiplicity}
Under the assumptions of Theorem~\ref{thm:existence},
\eqref{eq:original} has a sequence of radial finite-energy weak
solutions $\{u_j\}\subset E_V$ such that
\begin{equation*}
 0<J_0(u_1)<J_0(u_2)<\cdots\quad\text{and}\quad J_0(u_j)\longrightarrow+\infty\quad\text{as}\quad j\to\infty.
\end{equation*}
\end{theorem}

Two difficulties arise in the limit $\lambda\to0$. Each is handled by a
property of $f$.

The first difficulty is to bound critical points uniformly in $\lambda$.
A direct computation gives
\begin{equation}\label{eq:excess-intro}
 \frac12f(t)t-F(t)
 =\frac\alpha2\int_0^{|t|}s|\log s|^{\alpha-1}\,\mathrm{d}s
 \ge d_\alpha t^2,\qquad d_\alpha>0.
\end{equation}
In $J_\lambda(u)-\frac12\pair{J_\lambda'(u)}{u}$ all quadratic terms
cancel, including the potential term. So, for critical points,
\eqref{eq:excess-intro} bounds $|u|_2$ uniformly in $\lambda$. The
Gagliardo--Nirenberg inequality \cite{Nirenberg} then bounds
$|\nabla u|_2$. The same argument bounds Palais--Smale sequences of
$J_\lambda$ for fixed $\lambda$. No logarithmic Sobolev inequality is
needed.

The second difficulty is compactness. Since $V$ need not be coercive,
the energy space does not embed compactly into $L^2(\R^N)$ in general.
Here the infinite-mass structure helps. We have $-F(t)\ge\frac12
t^2|\log|t||^\alpha$ for $0<|t|\le1$. By the radial lemma of Strauss
\cite{Strauss}, radial functions in $H^1(\R^N)$ are small for large
$|x|$. Together these facts give
\begin{equation*}
 \sup_n\int_{|x|>R}|u_n|^2\dx
       \le \frac{C}{(\log R)^\alpha}\longrightarrow0
       \quad(R\to\infty)
\end{equation*}
for bounded sequences in $H^1_{\rad}(\R^N)$ whose energies are bounded
from above. By Lemma 2.6, if \(\lambda_n\to0^+\) and \(u_n\) are critical points of \(J_{\lambda_n}\) with energies bounded above uniformly in \(n\), then a subsequence converges strongly in \(E_V\) to a solution \(u\) of \eqref{eq:original}, and \(J_{\lambda_n}(u_n)\to J_{0}(u)\). The same perturbation was
used in \cite{HY} for $\alpha=1$. There the logarithmic Sobolev
inequality and the coercivity of the potential play an essential role.

\begin{remark}\label{rem:alpha}
The proofs of Theorems~\ref{thm:existence} and~\ref{thm:multiplicity}
also work for $\alpha=1$, where $\frac12f(t)t-F(t)=\frac14t^2$. We
assume $\alpha>1$ for the following reason. Let $\alpha=1$ and $\mu>0$.
Then $v$ solves $-\Delta v+V(x)v=v\log|v|$ if and only if $u=\mu^{-1}v$
solves
\begin{equation*}
 -\Delta u+(V(x)-\log\mu)u=u\log|u|.
\end{equation*}
So for $\alpha=1$ a potential bounded from below can be shifted to a
positive one. For $\alpha>1$ this scaling argument fails, since
$\mu^{-1}f(\mu t)-f(t)$ is not linear in $t$ when $\mu\ne1$. Thus, for
$\alpha>1$, sign-changing potentials cannot be reduced to positive ones
in this way.
\end{remark}

We next consider the more general problem
\begin{equation}\label{eq:general}
 -\Delta u+V(x)u=g(u)\quad\text{in }\R^N.
\end{equation}
The condition $g(t)/t\to-\infty$ as $t\to0$ is the infinite-mass case of
nonlinear scalar field equations. For the autonomous problem
$-\Delta u=g(u)$, Mederski \cite{Mederski} developed a general theory.
It gives least energy solutions, multiple nonradial solutions and a
sequence of radial solutions with unbounded energies. In particular, it
covers \eqref{eq:original} with a constant potential. Here we allow nonconstant radial potentials.

Let $S>0$ be the best constant in the Sobolev inequality
\begin{equation}\label{eq:Sobolev}
 S\left(\int_{\R^N}|u|^{2^*}\dx\right)^{2/2^*}
       \le\int_{\R^N}|\nabla u|^2\dx,
       \quad\text{for all } u\in C_c^\infty(\R^N).
\end{equation}
For \eqref{eq:general} we assume
\begin{equation}\tag{V2}\label{V2}
 \begin{gathered}
 V\in C^1(\R^N)\cap L^\infty(\R^N),\qquad V(x)=V(|x|),\\
 W\in L^\infty(\R^N),\qquad
 W_+\in L^{N/2}(\R^N),\qquad |W_+|_{N/2}<2S,
 \end{gathered}
\end{equation}
where $W(x):=x\cdot\nabla V(x)$. Let
\begin{equation*}
 G(t)=\int_0^t g(s)\,\mathrm{d}s,\qquad
 G_\pm(t)=\int_0^{|t|}(\pm g(s))_+\,\mathrm{d}s.
\end{equation*}
If $g$ is odd, then $G=G_+-G_-$ and $G_\pm\ge0$. Note that $G_+$ is in
general not the positive part of $G$. We assume:
\begin{equation}\tag{g1}\label{g1}
 g\in C(\R,\R),\quad g\text{ is odd},\quad
 \lim_{t\to0}\frac{g(t)}t=-\infty;
\end{equation}
\begin{equation}\tag{g2}\label{g2}
 \text{there are }r\in(1,2)\text{ and }C>0\text{ such that }
       |g(t)|\le C|t|^{r-1}\quad(|t|\le1);
\end{equation}
\begin{equation}\tag{g3}\label{g3}
 \limsup_{|t|\to\infty}\frac{|g(t)|}{|t|^{2^*-1}}<\infty,
       \qquad
 \lim_{|t|\to\infty}\frac{G_+(t)}{|t|^{2^*}}=0;
\end{equation}
\begin{equation}\tag{g4}\label{g4}
 \text{there is }\xi>0\text{ such that }
       G(\xi)>\frac12\Bigl(\sup_{\R^N}V\Bigr)\xi^2.
\end{equation}

\begin{remark}
By \eqref{g2}, once the term $\frac\lambda p\int_{\R^N}|u|^p\dx$ with
$1<p<r$ is added, the energy functional is finite and of class $C^1$ on
$H^1_{\rad}(\R^N)\cap L^p(\R^N)$. Condition \eqref{g2} is an extra
assumption compared with the hypotheses of \cite{Mederski} for the
autonomous problem. Condition \eqref{g3} allows critical growth. For
instance, a negative critical power term is allowed (see
Remark~\ref{rem:g-examples}). Condition \eqref{g4} is used to build odd
maps on which the energy is negative in the presence of the potential.
\end{remark}

We say that $u\in H^1_{\rad}(\R^N)$ is a radial finite-energy weak
solution of \eqref{eq:general} if $G(u)\in L^1(\R^N)$ and
\begin{equation}\label{eq:general-weak}
 \int_{\R^N}\bigl(\nabla u\cdot\nabla\varphi+V(x)u\varphi\bigr)\dx
       =\int_{\R^N}g(u)\varphi\dx,
       \quad\text{for all}\ \varphi\in C_c^\infty(\R^N).
\end{equation}
Its energy is
\begin{equation*}
 I_0(u)=\frac12\int_{\R^N}\bigl(|\nabla u|^2+V(x)u^2\bigr)\dx
                  -\int_{\R^N}G(u)\dx.
\end{equation*}

\begin{theorem}\label{thm:general-multiplicity}
Assume $N\ge3$, \eqref{V2} and \eqref{g1}--\eqref{g4}. Then
\eqref{eq:general} has a sequence of radial finite-energy weak
solutions $\{u_j\}$ such that
\begin{equation*}
 0<I_0(u_1)<I_0(u_2)<\cdots\quad\text{and}\quad I_0(u_j)\longrightarrow+\infty\quad\text{as}\quad j\to\infty.
\end{equation*}
Moreover, $g(u_j)u_j\in L^1(\R^N)$ for every $j$.
\end{theorem}

The proof of Theorem~\ref{thm:general-multiplicity} relies on the
following convergence result. Fix $1<p<r$ and set
\begin{equation}\label{eq:I-intro}
 I_\lambda(u)=I_0(u)+\frac\lambda p\int_{\R^N}|u|^p\dx,
       \qquad u\in H^1_{\rad}(\R^N)\cap L^p(\R^N).
\end{equation}

\begin{theorem}\label{thm:general-limit}
Assume $N\ge3$, \eqref{V2} and \eqref{g1}--\eqref{g3}. Let
$\lambda_n\to0^+$, and let $u_n\in H^1_{\rad}(\R^N)\cap L^p(\R^N)$ be
critical points of $I_{\lambda_n}$ with
$\sup_n I_{\lambda_n}(u_n)<\infty$. Then, after passing to a
subsequence, there is a radial finite-energy weak solution $u$ of
\eqref{eq:general} such that
\begin{equation*}
 \begin{gathered}
 u_n\to u\text{ in }H^1(\R^N),\qquad
 \lambda_n\int_{\R^N}|u_n|^p\dx\longrightarrow0,\\
 G(u_n)\to G(u)\text{ in }L^1(\R^N),\qquad
 I_{\lambda_n}(u_n)\longrightarrow I_0(u).
 \end{gathered}
\end{equation*}
If $\liminf_n I_{\lambda_n}(u_n)>0$, then $u\ne0$.
\end{theorem}

For general $g$, the lower bound in \eqref{eq:excess-intro} may fail
(see Remark~\ref{rem:g-examples}). Then it is not clear that
Palais--Smale sequences are bounded. Instead we use the Pohozaev
identity and the augmented functional method of Hirata, Ikoma and
Tanaka \cite{HIT}. This method uses the dilations
$u\mapsto u(\mathrm{e}^{-\theta}x)$. The Pohozaev identity gives bounds
that do not depend on $\lambda$ (Lemma~\ref{lem:general-bounds}). So we
can combine the perturbation $\lambda|u|^{p-2}u$ with the augmented
minimax method, and remove it in the limit $\lambda\to0$. No
Ambrosetti--Rabinowitz type condition is needed. The smallness condition
$|W_+|_{N/2}<2S$ was used by Azzollini and
Pomponio~\cite{AzzolliniPomponio2009} to obtain radial solutions under
Berestycki--Lions type assumptions. Theorem~\ref{thm:general-multiplicity}
extends the radial multiplicity result of \cite{Mederski} to nonconstant
radial potentials satisfying \eqref{V2}, under the extra condition
\eqref{g2}.

\begin{remark}\label{rem:potential-examples}
(i) Condition \eqref{V2} holds for
\begin{equation*}
 V(x)=c+\frac{a}{1+|x|^2},\qquad c\in\R,\quad a\ge0,
\end{equation*}
since $W(x)=-2a|x|^2(1+|x|^2)^{-2}$ is bounded and nonpositive. These
potentials change sign when $c<0<c+a$.

(ii) Let $V$ be as in (i) and $\varphi\in C_{c,\rad}^\infty(\R^N)$.
Then $V+\varepsilon\varphi$ satisfies \eqref{V2} whenever
$|\varepsilon|\,|x\cdot\nabla\varphi|_{N/2}<2S$. So nonmonotone
potentials are also allowed.
\end{remark}

\begin{remark}\label{rem:thm1andthm2}
The nonlinearity \eqref{eq:fF} satisfies \eqref{g1}--\eqref{g4} for
every bounded $V$. So Theorem~\ref{thm:general-multiplicity} also
applies to \eqref{eq:original} when $V$ satisfies \eqref{V2}. For
$u\in L^2(\R^N)$,
\begin{equation*}
 F(u)\in L^1(\R^N)\iff\int_{\R^N}|u|^2|\log|u||^\alpha\dx<\infty,
\end{equation*}
by \eqref{eq:F-domination}, \eqref{eq:F-small} and
\eqref{eq:F-asymptotic}. Hence the two notions of radial finite-energy
weak solution coincide for \eqref{eq:original}.
Theorems~\ref{thm:existence} and~\ref{thm:multiplicity}, however, hold
for a much larger class of potentials than \eqref{V2}.
\end{remark}

\begin{remark}\label{rem:g-examples}
Theorem~\ref{thm:general-multiplicity} applies, for instance, to
\begin{equation*}
 g(t)=-|t|^{r-2}t+a|t|^{q-2}t-|t|^{2^*-2}t,
       \qquad 1<r<2<q<2^*,
\end{equation*}
provided
\begin{equation*}
 a>\max\left\{0,\ q\left(\frac1r+\frac1{2^*}
              +\frac12\sup_{\R^N}V\right)\right\}.
\end{equation*}
For this $g$,
\begin{equation*}
 \frac12g(t)t-G(t)
 =\left(\frac1r-\frac12\right)|t|^r
  +a\left(\frac12-\frac1q\right)|t|^q
  -\left(\frac12-\frac1{2^*}\right)|t|^{2^*}
       \longrightarrow-\infty\quad(|t|\to\infty).
\end{equation*}
So the lower bound in \eqref{eq:excess-intro} fails, and the arguments
of Sections~\ref{sec:existence} and~\ref{sec:multiplicity} do not apply
directly.
\end{remark}

The paper is organized as follows. In Section~\ref{sec:existence} we
prove Theorem~\ref{thm:existence}. We study $f$, verify the
Palais--Smale condition for $J_\lambda$ and pass to the limit
$\lambda\to0$. In Section~\ref{sec:multiplicity} we prove
Theorem~\ref{thm:multiplicity}. In Section~\ref{sec:general} we treat
general nonlinearities and prove Theorems~\ref{thm:general-multiplicity}
and~\ref{thm:general-limit}.

\section{Proof of Theorem \ref{thm:existence}}\label{sec:existence}

Throughout this section we assume \eqref{V1}. Let $X$ and $J_\lambda$ be
given by \eqref{eq:X} and \eqref{eq:energy}, with $1<p<2$ fixed. The map
$u\mapsto(u,u)$ is an isometry from $X$ onto a closed subspace of
$E_V\times L^p(\R^N)$ with the norm $\norm{(v,w)}=\norm{v}_V+|w|_p$. So
$X$ is separable and reflexive. The function $V_+$ is bounded on bounded
sets. So radial cut-off functions and radial mollifiers show that
$C_{c,\rad}^\infty(\R^N)$ is dense in $X$ and in $E_V$. By interpolation
and the Sobolev inequality, $X\hookrightarrow L^s(\R^N)$ for
$p\le s\le2^*$.

\begin{lemma}\label{lem:fF}
Let \(f\) be defined by \eqref{eq:fF} and \(F(t)=\int_0^t f(s)\,\mathrm{d}s\).
Then \(f\) is continuous and odd, \(F\in C^1(\R)\) is even, and the
following hold.

\begin{enumerate}
\item[(i)] For every sufficiently small \(\varepsilon>0\), there
is \(C_\varepsilon>0\) such that, for all \(t\in\R\),
\begin{equation}\label{eq:growth}
 \begin{split}
 |f(t)|&\le C_\varepsilon
        (|t|^{1-\varepsilon}+|t|^{1+\varepsilon}),\\
 |F(t)|+|f(t)t|&\le C_\varepsilon
        (|t|^{2-\varepsilon}+|t|^{2+\varepsilon}).
 \end{split}
\end{equation}

\item[(ii)] For all \(t\in\R\),
\begin{equation}\label{eq:excess}
 \frac12 f(t)t-F(t)
       =\frac{\alpha}{2}\int_0^{|t|}s|\log s|^{\alpha-1}\,\mathrm{d}s.
\end{equation}
There are \(d_\alpha>0\) and \(C>0\) such that
\begin{equation}\label{eq:coercive}
 \frac12f(t)t-F(t)\ge d_\alpha t^2
\end{equation}
and
\begin{equation}\label{eq:F-domination}
 |F(t)|+\left(\frac12f(t)t-F(t)\right)
       \le C\bigl(t^2+|f(t)t|\bigr).
\end{equation}

\item[(iii)] For \(0<|t|\le1\),
\begin{equation}\label{eq:F-small}
 -F(t)=\int_0^{|t|}s|\log s|^\alpha\,\mathrm{d}s
       \ge\frac12t^2|\log|t||^\alpha.
\end{equation}
For every \(A>0\) and \(q>2\), there are \(C_{A,q}>0\) and \(C_q>0\)
such that, for all \(t\in\R\),
\begin{equation}\label{eq:F-upper}
 F(t)\le-At^2+C_{A,q}|t|^q,
       \qquad (f(t)t)_+\le C_q|t|^q.
\end{equation}

\item[(iv)]
\begin{equation}\label{eq:F-asymptotic}
 \lim_{t\to\infty}\frac{F(t)}{t^2(\log t)^\alpha}=\frac12,
       \qquad\lim_{t\to0}\frac{F(t)}{|t|^p}=0.
\end{equation}
For every \(L>0\), there is \(C_L>0\) such that
\begin{equation}\label{eq:F-lower}
 F(t)\ge Lt^2-C_L|t|^p,\quad\text{for all}~~t\in\R.
\end{equation}
\end{enumerate}
\end{lemma}

\begin{proof}
Since \(|t||\log|t||^\alpha\to0\) as \(t\to0\), \(f\) is continuous.
It is odd. So \(F\) is even and \(F\in C^1(\R)\).

(i) Let \(\varepsilon>0\). The function \(s^\varepsilon|\log s|^\alpha\)
is bounded on \((0,1]\), and \(s^{-\varepsilon}(\log s)^\alpha\) is
bounded on \([1,\infty)\). This gives the bound for \(f\). We multiply
by \(|t|\) and integrate. This gives the bounds for \(f(t)t\) and \(F\).

(ii) Let \(t>0\). We have
\(\frac{\mathrm d}{\mathrm ds}\bigl(|\log s|^{\alpha-1}\log s\bigr)
=\alpha s^{-1}|\log s|^{\alpha-1}\) and
\(s^2|\log s|^\alpha\to0\) as \(s\to0^{+}\). So integration by parts
gives
\begin{equation*}
 F(t)=\frac12f(t)t-\frac\alpha2\int_0^ts|\log s|^{\alpha-1}\,\mathrm{d}s,
\end{equation*}
which is \eqref{eq:excess}. By evenness, \eqref{eq:excess} holds for all
\(t\). Its right-hand side is positive for \(t\ne0\). If \(0<|t|\le1\),
then \(|\log s|\ge|\log|t||\) for \(0<s\le|t|\). Hence
\begin{equation}\label{eq:excess-small}
 \frac12f(t)t-F(t)
       \ge\frac{\alpha}{4}t^2|\log|t||^{\alpha-1},
       \qquad 0<|t|\le1.
\end{equation}
For \(t\ge4\), we integrate only over \([t/2,t]\) and get
\begin{equation*}
 \frac{\frac12f(t)t-F(t)}{t^2}
       \ge\frac{\alpha}{2t^2}\Bigl(\log\frac t2\Bigr)^{\alpha-1}
          \int_{t/2}^ts\,\mathrm{d}s
       =\frac{3\alpha}{16}\Bigl(\log\frac t2\Bigr)^{\alpha-1}
       \ge\frac{3\alpha}{16}(\log2)^{\alpha-1}.
\end{equation*}
By \eqref{eq:excess-small}, the same quotient is at least
\(\frac\alpha4\) for \(0<t\le\mathrm{e}^{-1}\). The quotient is
continuous and positive on \((0,\infty)\). So its infimum is positive.
By evenness, \eqref{eq:coercive} follows.

For \(0<|t|\le1\), the substitution \(s=|t|r\) gives
\begin{equation}\label{eq:excess-upper-proof}
 \begin{split}
 \frac12f(t)t-F(t)
 &=\frac{\alpha t^2}{2}\int_0^1
       r\bigl(|\log|t||+|\log r|\bigr)^{\alpha-1}\,\mathrm{d}r\\
 &\le Ct^2\bigl(1+|\log|t||^{\alpha-1}\bigr),
 \end{split}
\end{equation}
since \(\int_0^1r|\log r|^{\alpha-1}\,\mathrm{d}r<\infty\). For
\(|t|\ge1\), we split the integral in \eqref{eq:excess} at \(s=1\). This
gives the same bound. Now \eqref{eq:F-domination} follows from
\(s^{\alpha-1}\le C(1+s^\alpha)\), \(|f(t)t|=t^2|\log|t||^\alpha\) and
\[
 F(t)=\tfrac12f(t)t-\bigl(\tfrac12f(t)t-F(t)\bigr).
\]

(iii) On \((0,1)\), \(f(s)=-s|\log s|^\alpha\). Also
\(|\log s|\ge|\log|t||\) for \(0<s\le|t|\le1\). This gives
\eqref{eq:F-small}. In particular, $F(t)/t^2\to-\infty$ as $t\to0$. Fix
$A>0$ and $q>2$. Then $(F(t)+At^2)_+=0$ near $t=0$. By \eqref{eq:growth}
with $\varepsilon<q-2$,
\[
 \frac{(F(t)+At^2)_+}{|t|^q}\longrightarrow0\quad(|t|\to\infty).
\]
So this continuous function on \(\{t\ne0\}\) is bounded. This is the
first inequality in \eqref{eq:F-upper}. The second one follows from
\(f(t)t\le0\) on \([-1,1]\) and \((\log t)^\alpha\le C_qt^{q-2}\) for
\(t\ge1\).

(iv) The first limit in \eqref{eq:F-asymptotic} follows from
l'H\^opital's rule. The second one follows from \eqref{eq:growth} with
\(\varepsilon<2-p\). Let \(L>0\). The function
\((Lt^2-F(t))_+/|t|^p\), \(t\ne0\), is continuous. It tends to zero as
\(t\to0\), and it vanishes for \(|t|\) large. So it is bounded, and
\eqref{eq:F-lower} follows.
\end{proof}

By \eqref{eq:F-domination}, $J_0(u)$ is finite for every $u\in E_V$ with
$\int_{\R^N}|u|^2|\log|u||^\alpha\dx<\infty$. In particular, a weak
solution in the sense of \eqref{eq:weak} has finite energy.

We show that \(J_\lambda\in C^1(X,\R)\). Choose \(\varepsilon>0\) with
\(p<2(1-\varepsilon)<2<2(1+\varepsilon)<2^*\). By \eqref{eq:growth},
\(J_\lambda\) is finite on \(X\), and
\begin{equation}\label{eq:f-L2}
 |f(u)|_2\le C\bigl(|u|_{2(1-\varepsilon)}^{1-\varepsilon}
                   +|u|_{2(1+\varepsilon)}^{1+\varepsilon}\bigr).
\end{equation}
Let \(u_n\to u\) in \(X\). Then
\(|u_n|^{2(1\pm\varepsilon)}\to|u|^{2(1\pm\varepsilon)}\) in
\(L^1(\R^N)\). We apply Vitali's theorem along subsequences that
converge almost everywhere. This gives \(f(u_n)\to f(u)\) in
\(L^2(\R^N)\). For \(u,h\in X\), the fundamental theorem of calculus and
Fubini's theorem give
\begin{equation}\label{eq:differentiability}
 \begin{split}
 \left|\int_{\R^N}\bigl(F(u+h)-F(u)-f(u)h\bigr)\dx\right|
 &\le |h|_2\sup_{0\le s\le1}|f(u+sh)-f(u)|_2\\
 &=o(\norm{h}_X).
 \end{split}
\end{equation}
The same argument shows that the derivative is continuous. The potential
term is continuous on $E_V$. This follows from the Cauchy--Schwarz
inequality with weight $V_+$ and from $0\le V_-\le M_0$. Thus
\(J_\lambda\in C^1(X,\R)\), \(J_\lambda\) is even, and
\begin{equation*}
 \pair{J_\lambda'(u)}{v}
 =\int_{\R^N}(\nabla u\cdot\nabla v+V(x)uv)\dx
        +\lambda\int_{\R^N}|u|^{p-2}uv\dx-\int_{\R^N}f(u)v\dx,
 \qquad u,v\in X.
\end{equation*}

Next we bound the positive part of the nonlinear term. Let
\begin{equation*}
 \delta=\frac{4(p-1)}{N(p-1)+4}.
\end{equation*}
Then \(0<\delta<\min\{4/N,2^*-2\}\) and
\begin{equation*}
 \frac{2+\delta-N\delta/2}{1-N\delta/4}=p+1.
\end{equation*}
By \eqref{eq:F-upper} with \(A=1\) and \(q=2+\delta\),
\(F(t)\le C|t|^{2+\delta}\). We use the Gagliardo--Nirenberg inequality
\cite{Nirenberg} and then Young's inequality with exponents
\(4/(N\delta)\) and \(1/(1-N\delta/4)\). This gives
\begin{equation}\label{eq:GN-F}
 \int_{\R^N}F(u)\dx
 \le C\int_{\R^N}|u|^{2+\delta}\dx
 \le C|\nabla u|_2^{N\delta/2}|u|_2^{2+\delta-N\delta/2}
 \le\frac14|\nabla u|_2^2+C|u|_2^{p+1},
 \qquad u\in X.
\end{equation}

Recall that the embedding \(H^1_{\rad}(\R^N)\hookrightarrow L^q(\R^N)\)
is compact for \(2<q<2^*\) \cite{Strauss}.

\begin{lemma}\label{lem:compactness}
Suppose that \(u_n\in X\), \(0\le\lambda_n\le1\), and
\begin{equation}\label{eq:compactness-assumption}
 \sup_n\norm{u_n}<\infty,\qquad
       \sup_nJ_{\lambda_n}(u_n)<\infty.
\end{equation}
Then, for all sufficiently large \(R\),
\begin{equation}\label{eq:log-tail}
 \sup_n\int_{|x|>R}|u_n|^2\dx
       \le\frac{C}{(\log R)^\alpha}.
\end{equation}
Moreover, \(\{u_n\}\) is relatively compact in \(L^2(\R^N)\).
\end{lemma}

\begin{proof}
Radial functions are small for large $|x|$. Near
$t=0$, $-F(t)$ is much larger than $t^2$. So a bound on the energy
controls the $L^2$ norm of $u_n$ outside a large ball.

Let \(u(x)=U(|x|)\) be radial, smooth and compactly supported. Then
\begin{equation*}
 r^{N-1}U(r)^2
 =-\int_r^\infty\frac{\mathrm{d}}{\mathrm{d}s}(s^{N-1}U(s)^2)\,\mathrm{d}s
 \le2\int_r^\infty s^{N-1}|U(s)U'(s)|\,\mathrm{d}s
 \le C_N|u|_2|\nabla u|_2.
\end{equation*}
By density, every radial \(u\in H^1(\R^N)\) has a representative with
\begin{equation}\label{eq:radial}
 |u(x)|\le C_N|x|^{-(N-1)/2}\norm{u},
       \quad\text{for almost every }x\ne0.
\end{equation}

Fix \(q\in(2,2^*)\). By \eqref{eq:F-small}, \(F(t)\le0\) for
\(|t|\le1\). By \eqref{eq:F-upper} with \(A=1\),
\(F(t)\le C_{1,q}|t|^q\). We split \(\int_{\R^N}F(u_n)\dx\) into
integrals over \(\{|u_n|\le1\}\) and \(\{|u_n|>1\}\). Then
\eqref{eq:energy0}, \eqref{eq:energy} and the Sobolev inequality give
\begin{equation}\label{eq:negative-energy-bound}
 \begin{split}
 &\int_{\{|u_n|\le1\}}\bigl(-F(u_n)\bigr)\dx
       +\frac12\int_{\R^N}V_+(x)u_n^2\dx
       +\frac{\lambda_n}{p}|u_n|_p^p\\
 &\quad=J_{\lambda_n}(u_n)-\frac12|\nabla u_n|_2^2
       +\frac12\int_{\R^N}V_-(x)u_n^2\dx
       +\int_{\{|u_n|>1\}}F(u_n)\dx\\
 &\quad\le C+\frac{M_0}{2}|u_n|_2^2+C_{1,q}|u_n|_q^q\le C.
 \end{split}
\end{equation}
All three terms on the left are nonnegative. By \eqref{eq:radial}, for
\(R\) large and \(|x|>R\),
\begin{equation}\label{eq:log-lower-tail}
 |u_n(x)|\le \mathrm{e}^{-1},\qquad
 |\log|u_n(x)||\ge\frac{N-1}{4}\log R
       \quad\text{if }u_n(x)\ne0.
\end{equation}
With \eqref{eq:F-small}, this gives
\begin{equation*}
 -F(u_n(x))\ge
       \frac12\left(\frac{N-1}{4}\log R\right)^\alpha
       |u_n(x)|^2,\qquad |x|>R.
\end{equation*}
Since \(\{|x|>R\}\subset\{|u_n|\le1\}\) up to a null set, we integrate
over \(\{|x|>R\}\) and use \eqref{eq:negative-energy-bound}. This gives
\eqref{eq:log-tail}.

Every subsequence of \(\{u_n\}\) has a further subsequence which
converges to some \(u\) weakly in \(H^1(\R^N)\), strongly in
\(L^2_{\mathrm{loc}}(\R^N)\) and almost everywhere. By Fatou's lemma,
\(u\) also satisfies \eqref{eq:log-tail}. We split \(|u_n-u|_2^2\) into
integrals over \(B_R\) and \(\R^N\setminus B_R\). This gives
\(u_n\to u\) in \(L^2(\R^N)\).
\end{proof}

\begin{lemma}\label{lem:PS}
For every fixed $0<\lambda\le1$, $J_\lambda$ satisfies the
Palais--Smale condition on $X$.
\end{lemma}

\begin{proof}
Let $|J_\lambda(u_n)|\le C$ and $J_\lambda'(u_n)\to0$ in $X^*$. First we
show that $\{u_n\}$ is bounded. For $u\in X$,
\begin{equation}\label{eq:energy-derivative}
 J_\lambda(u)-\frac12\pair{J_\lambda'(u)}{u}
 =\int_{\R^N}\left(\frac12f(u)u-F(u)\right)\dx
       +\lambda\left(\frac1p-\frac12\right)|u|_p^p.
\end{equation}
By \eqref{eq:coercive},
\begin{equation*}
 d_\alpha|u_n|_2^2
       +\lambda\left(\frac1p-\frac12\right)|u_n|_p^p
       \le C+o(1)\norm{u_n}_X.
\end{equation*}
Since $\lambda>0$ is fixed,
\begin{equation}\label{eq:PS-L2}
 |u_n|_2^2+|u_n|_p^p\le C_\lambda(1+\norm{u_n}_X).
\end{equation}
By \eqref{eq:energy0}, \eqref{eq:energy}, \eqref{eq:GN-F} and
$V_-\le M_0$,
\begin{equation*}
 \frac14|\nabla u_n|_2^2+\frac12\int_{\R^N}V_+(x)u_n^2\dx
 \le J_\lambda(u_n)+\frac{M_0}{2}|u_n|_2^2+C|u_n|_2^{p+1}.
\end{equation*}
By \eqref{eq:PS-L2}, the right-hand side is at most
$C_\lambda(1+\norm{u_n}_X^{(p+1)/2})$, since $(p+1)/2>1$. Hence
\begin{equation*}
 \norm{u_n}_V^2
 \le4\left(\frac14|\nabla u_n|_2^2
           +\frac12\int_{\R^N}V_+(x)u_n^2\dx\right)+|u_n|_2^2
 \le C_\lambda\bigl(1+\norm{u_n}_X^{(p+1)/2}\bigr).
\end{equation*}
Since $1<p<2$, it follows that
\begin{equation*}
 \norm{u_n}_X^p
 \le C\bigl(1+\norm{u_n}_V^2+|u_n|_p^p\bigr)
 \le C_\lambda\bigl(1+\norm{u_n}_X^{(p+1)/2}\bigr).
\end{equation*}
As $(p+1)/2<p$, $\{u_n\}$ is bounded in $X$.

By Lemma~\ref{lem:compactness}, up to a subsequence, $u_n\weakto u$ in
$X$ and $u_n\to u$ in $L^2(\R^N)$. By \eqref{eq:f-L2}, $\{f(u_n)\}$ is
bounded in $L^2(\R^N)$. So
\begin{equation*}
 \begin{gathered}
 \left|\int_{\R^N}f(u_n)(u_n-u)\dx\right|
       \le |f(u_n)|_2|u_n-u|_2\to0,\\
 \left|\int_{\R^N}V_-(x)u_n(u_n-u)\dx\right|
       \le M_0|u_n|_2|u_n-u|_2\to0.
 \end{gathered}
\end{equation*}
We test $J_\lambda'(u_n)$ with $u_n-u$. Then we subtract
$\int_{\R^N}(\nabla u\cdot\nabla(u_n-u)+V_+(x)u(u_n-u))\dx$ and
$\lambda\int_{\R^N}|u|^{p-2}u(u_n-u)\dx$, which tend to zero by weak
convergence. Since also $|u_n-u|_2\to0$, we get
\begin{equation}\label{eq:PS-monotone}
 \norm{u_n-u}_V^2
 +\lambda\int_{\R^N}
  \bigl(|u_n|^{p-2}u_n-|u|^{p-2}u\bigr)(u_n-u)\dx=o(1).
\end{equation}
Both terms on the left are nonnegative. For $a,b\in\R$ with
$|a|+|b|>0$,
\begin{equation*}
 \bigl(|a|^{p-2}a-|b|^{p-2}b\bigr)(a-b)
 \ge(p-1)\frac{|a-b|^2}{(|a|+|b|)^{2-p}}.
\end{equation*}
Hence, by H\"older's inequality with exponents $2/p$ and $2/(2-p)$,
\begin{equation*}
 \begin{split}
 |u_n-u|_p^p
 &\le\left(\int_{\R^N}
       \frac{|u_n-u|^2}{(|u_n|+|u|)^{2-p}}\dx\right)^{p/2}
       \left(\int_{\R^N}(|u_n|+|u|)^p\dx\right)^{(2-p)/2}\\
 &\le C\left(\int_{\R^N}
    (|u_n|^{p-2}u_n-|u|^{p-2}u)(u_n-u)\dx\right)^{p/2}
       \longrightarrow0.
 \end{split}
\end{equation*}
With \eqref{eq:PS-monotone}, this gives $u_n\to u$ in $X$.
\end{proof}

\begin{lemma}\label{lem:geometry}
There exist \(\rho>0\), \(m_1>0\) and \(e_0\in X\), independent of
\(0<\lambda\le1\), such that
\begin{equation*}
 \inf_{\substack{u\in X\\\norm{u}_V=\rho}}J_\lambda(u)\ge m_1,
       \qquad \norm{e_0}_V>\rho,\qquad J_\lambda(e_0)<0.
\end{equation*}
\end{lemma}

\begin{proof}
Fix \(q\in(2,2^*)\). We use \eqref{eq:F-upper} with $A=M_0/2+1/4$. Then
\begin{equation*}
 \frac12\int_{\R^N}V(x)u^2\dx-\int_{\R^N}F(u)\dx
 \ge\frac12\int_{\R^N}V_+(x)u^2\dx+\frac14|u|_2^2-C|u|_q^q.
\end{equation*}
With the Sobolev embedding, there is $C_0>0$ such that
\begin{equation}\label{eq:lower-geometry}
 \begin{split}
 J_\lambda(u)
 &\ge\frac14\norm{u}_V^2+\frac\lambda p|u|_p^p-C_0|u|_q^q\\
 &\ge\frac14\norm{u}_V^2+\frac\lambda p|u|_p^p-C\norm{u}_V^q.
 \end{split}
\end{equation}
Choose \(\rho>0\) with \(C\rho^{q-2}\le1/8\). Then
\begin{equation}\label{eq:small-ball}
 J_\lambda(u)\ge\frac18\norm{u}_V^2+\frac\lambda p|u|_p^p
       \qquad\text{if }\norm{u}_V\le\rho.
\end{equation}
So the first assertion holds with \(m_1=\rho^2/8\).

Take \(\varphi\in C_{c,\rad}^\infty(\R^N)\setminus\{0\}\). By
\eqref{eq:F-lower} and $\lambda\le1$, for every \(L>0\),
\begin{equation*}
 J_\lambda(t\varphi)
 \le t^2\left(\frac12\int_{\R^N}(|\nabla\varphi|^2+V(x)\varphi^2)\dx-L|\varphi|_2^2\right)
       +\left(C_L+\frac1p\right)t^p|\varphi|_p^p.
\end{equation*}
Choose \(L\) such that the coefficient of \(t^2\) is negative. Since
\(p<2\), the right-hand side tends to \(-\infty\) as \(t\to\infty\),
uniformly in \(\lambda\in(0,1]\). Take \(e_0=t_0\varphi\) with \(t_0\)
large.
\end{proof}

Let
\begin{equation}\label{eq:MP-level}
 \Gamma=\{\gamma\in C([0,1],X):\gamma(0)=0,\ \gamma(1)=e_0\},
       \qquad c_\lambda=\inf_{\gamma\in\Gamma}
                         \max_{s\in[0,1]}J_\lambda(\gamma(s)).
\end{equation}
The mountain pass geometry also holds for the norm of \(X\). Indeed, let
\(\norm{u}_X=\rho\). Then \(\norm{u}_V\le\rho\), and either
\(\norm{u}_V\ge\rho/2\) or \(|u|_p\ge\rho/2\). By \eqref{eq:small-ball},
\begin{equation*}
 J_\lambda(u)\ge\min\left\{\frac{\rho^2}{32},\
       \frac\lambda p\Bigl(\frac\rho2\Bigr)^p\right\}>0.
\end{equation*}
Moreover, \(\norm{e_0}_X\ge\norm{e_0}_V>\rho\) and \(J_\lambda(0)=0\).
By the mountain pass theorem \cite{Willem}
and Lemma~\ref{lem:PS}, there is \(u_\lambda\in X\) with
\begin{equation}\label{eq:MP-critical}
 J_\lambda'(u_\lambda)=0,
       \qquad J_\lambda(u_\lambda)=c_\lambda.
\end{equation}
Every \(\gamma\in\Gamma\) crosses the sphere \(\norm{u}_V=\rho\), since
\(\norm{\gamma(0)}_V=0\) and \(\norm{\gamma(1)}_V>\rho\). So
Lemma~\ref{lem:geometry} gives \(c_\lambda\ge m_1\).

\begin{lemma}\label{lem:MP-bounds}
There is \(m_2>0\), independent of \(0<\lambda\le1\), such that
\begin{equation*}
 0<m_1\le J_\lambda(u_\lambda)\le m_2.
\end{equation*}
\end{lemma}

\begin{proof}
The lower bound was shown above. The path \(s\mapsto se_0\) belongs to
\(\Gamma\), and \(J_\lambda\le J_1\). Therefore
\begin{equation*}
 c_\lambda\le\max_{s\in[0,1]}J_\lambda(se_0)
       \le\max_{s\in[0,1]}J_1(se_0)=:m_2<\infty.
\end{equation*}
\end{proof}

The next lemma is the key step. It shows how to pass to the limit
$\lambda\to0$.

\begin{lemma}\label{lem:limit}
Let \(\lambda_n\in(0,1]\), \(\lambda_n\to0\), and let \(u_n\in X\)
satisfy, for some \(M>0\),
\begin{equation*}
 J_{\lambda_n}'(u_n)=0,\qquad J_{\lambda_n}(u_n)\le M.
\end{equation*}
Then, after passing to a subsequence, there is a radial finite-energy
weak solution \(u\in E_V\) of \eqref{eq:original} such that
\begin{equation}\label{eq:strong-limit}
 \begin{gathered}
 u_n\to u\text{ in }E_V,\qquad \lambda_n|u_n|_p^p\longrightarrow0,\\
 F(u_n)\to F(u)\text{ in }L^1(\R^N),\qquad
 J_{\lambda_n}(u_n)\longrightarrow J_0(u).
 \end{gathered}
\end{equation}
If \(\liminf_nJ_{\lambda_n}(u_n)>0\), then \(u\ne0\).
\end{lemma}

\begin{proof}
The proof has four steps.

\emph{Step 1.} By \eqref{eq:energy-derivative} and
\eqref{eq:coercive},
\begin{equation*}
 \begin{split}
 M\ge J_{\lambda_n}(u_n)
 &=\int_{\R^N}\left(\frac12f(u_n)u_n-F(u_n)\right)\dx
       +\lambda_n\left(\frac1p-\frac12\right)|u_n|_p^p\\
 &\ge d_\alpha|u_n|_2^2
       +\lambda_n\left(\frac1p-\frac12\right)|u_n|_p^p.
 \end{split}
\end{equation*}
So \(|u_n|_2^2+\lambda_n|u_n|_p^p\le C_M\). As in the proof of
Lemma~\ref{lem:PS}, \eqref{eq:GN-F} gives
\begin{equation*}
 \frac14|\nabla u_n|_2^2+\frac12\int_{\R^N}V_+(x)u_n^2\dx
       \le J_{\lambda_n}(u_n)+\frac{M_0}{2}|u_n|_2^2
                         +C|u_n|_2^{p+1}\le C_M.
\end{equation*}
Thus
\begin{equation}\label{eq:uniform-bounds}
 \norm{u_n}_V^2+\lambda_n|u_n|_p^p
       +\int_{\R^N}\left(\frac12f(u_n)u_n-F(u_n)\right)\dx\le C_M.
\end{equation}
By Lemma~\ref{lem:compactness} and interpolation, up to a subsequence,
\begin{equation}\label{eq:weak-limit}
 u_n\weakto u\text{ in }E_V,\qquad
 u_n\to u\text{ in }L^s(\R^N)\quad(2\le s<2^*),\qquad
 u_n\to u\text{ a.e. in }\R^N.
\end{equation}

\emph{Step 2.} Since \(J_{\lambda_n}'(u_n)=0\),
\begin{equation}\label{eq:local-perturbed}
 \int_{\R^N}(\nabla u_n\cdot\nabla\varphi+V(x)u_n\varphi)\dx
       +\lambda_n\int_{\R^N}|u_n|^{p-2}u_n\varphi\dx
       =\int_{\R^N}f(u_n)\varphi\dx
\end{equation}
for all \(\varphi\in C_{c,\rad}^\infty(\R^N)\). In fact
\eqref{eq:local-perturbed} holds for all
\(\varphi\in C_c^\infty(\R^N)\). Indeed, \(u_n\) and \(V\) are radial.
So no integral in \eqref{eq:local-perturbed} changes if \(\varphi\) is
replaced by \(\varphi(Q\,\cdot)\), \(Q\in O(N)\). Hence we may replace
\(\varphi\) by
\begin{equation*}
 \overline\varphi=\int_{O(N)}\varphi(Q\,\cdot)\,\mathrm{d}Q
 \in C_{c,\rad}^\infty(\R^N),
\end{equation*}
where \(\mathrm dQ\) is the normalized Haar measure. This is the principle
of symmetric criticality; see \cite{Willem}.

We pass to the limit in \eqref{eq:local-perturbed}. By H\"older's
inequality and \eqref{eq:uniform-bounds},
\begin{equation*}
 \left|\lambda_n\int_{\R^N}|u_n|^{p-2}u_n\varphi\dx\right|
 \le\lambda_n^{1/p}
       (\lambda_n|u_n|_p^p)^{(p-1)/p}|\varphi|_p
       \longrightarrow0.
\end{equation*}
Fix \(s\in(2,2^*)\). Then \(|f(t)|^2\le C_s(1+|t|^s)\). So
\(\{f(u_n)\}\) is bounded in \(L^2(B_R)\) for every \(R>0\). By Vitali's
theorem, \(f(u_n)\to f(u)\) in \(L^1_{\mathrm{loc}}(\R^N)\). Letting
\(n\to\infty\) in \eqref{eq:local-perturbed}, we obtain
\begin{equation}\label{eq:limit-local}
 \int_{\R^N}(\nabla u\cdot\nabla\varphi+V(x)u\varphi)\dx=\int_{\R^N}f(u)\varphi\dx,
       \quad\text{for all}~~\varphi\in C_c^\infty(\R^N).
\end{equation}

Next we test \(J_{\lambda_n}'(u_n)=0\) with \(u_n\). This gives
\begin{equation}\label{eq:critical-Nehari}
 \begin{split}
 &|\nabla u_n|_2^2+\int_{\R^N}V_+(x)u_n^2\dx
       +\lambda_n|u_n|_p^p+\int_{\R^N}(f(u_n)u_n)_-\dx\\
 &\qquad=\int_{\R^N}(f(u_n)u_n)_+\dx
                  +\int_{\R^N}V_-(x)u_n^2\dx.
 \end{split}
\end{equation}
By \eqref{eq:F-upper} and \eqref{eq:uniform-bounds}, the right-hand side
is bounded. We use Fatou's lemma for the negative part and
\eqref{eq:F-upper} for the positive part. This gives
\begin{equation}\label{eq:limit-integrability}
 \int_{\R^N}|f(u)u|\dx=\int_{\R^N}|u|^2|\log|u||^\alpha\dx<\infty.
\end{equation}
Thus \(u\) is a radial finite-energy weak solution of
\eqref{eq:original}.

\emph{Step 3.} We claim that
\begin{equation}\label{eq:limit-Nehari}
 \int_{\R^N}(|\nabla u|^2+V(x)u^2)\dx=\int_{\R^N}f(u)u\dx.
\end{equation}
Let \(\eta_R\in C_{c,\rad}^\infty(\R^N)\) satisfy \(0\le\eta_R\le1\),
\(\eta_R=1\) on \(B_R\), \(\eta_R=0\) outside \(B_{2R}\) and
\(|\nabla\eta_R|\le C/R\). Since \(f(u)\in L^2(B_{2R})\),
\eqref{eq:limit-local} extends by density to test functions in
\(H_0^1(B_{2R})\). Testing with \(\eta_Ru\), we get
\begin{equation*}
 \int_{\R^N}\eta_R(|\nabla u|^2+V(x)u^2)\dx+\int_{\R^N}u\nabla u\cdot\nabla\eta_R\dx
       =\int_{\R^N}\eta_R f(u)u\dx.
\end{equation*}
The second integral is bounded by \(C|u|_2|\nabla u|_2/R\). We let
\(R\to\infty\). By \eqref{eq:limit-integrability},
$\int_{\R^N}|V(x)|u^2\dx<\infty$ and dominated convergence, we obtain
\eqref{eq:limit-Nehari}.

\emph{Step 4.} We compare \eqref{eq:critical-Nehari}
with \eqref{eq:limit-Nehari}. Fix $q\in(2,2^*)$. Since
$(f(t)t)_+\le C_q|t|^q$, and $u_n\to u$ in $L^q(\R^N)$ and almost
everywhere, Vitali's theorem gives
\begin{equation}\label{eq:positive-part-convergence}
 (f(u_n)u_n)_+\longrightarrow(f(u)u)_+\quad\text{in }L^1(\R^N).
\end{equation}
Since $u_n\to u$ in $L^2(\R^N)$ and $0\le V_-\le M_0$,
\begin{equation}\label{eq:potential-negative-conv}
 \int_{\R^N}V_-(x)u_n^2\dx\longrightarrow
       \int_{\R^N}V_-(x)u^2\dx.
\end{equation}
So, by \eqref{eq:limit-Nehari}, the right-hand side of
\eqref{eq:critical-Nehari} tends to
\begin{equation*}
 \int_{\R^N}(f(u)u)_+\dx+\int_{\R^N}V_-(x)u^2\dx
 =|\nabla u|_2^2+\int_{\R^N}V_+(x)u^2\dx+\int_{\R^N}(f(u)u)_-\dx.
\end{equation*}
On the other hand, by weak lower semicontinuity and Fatou's lemma,
\begin{equation*}
 \begin{gathered}
 |\nabla u|_2^2\le\liminf_{n\to\infty}|\nabla u_n|_2^2,\qquad
 \int_{\R^N}V_+u^2\dx\le\liminf_{n\to\infty}\int_{\R^N}V_+u_n^2\dx,\\
 \int_{\R^N}(f(u)u)_-\dx\le\liminf_{n\to\infty}\int_{\R^N}(f(u_n)u_n)_-\dx,
 \end{gathered}
\end{equation*}
and $\lambda_n|u_n|_p^p\ge0$. So the sum of the four terms on the left
of \eqref{eq:critical-Nehari} converges to the sum of their lower
bounds. Therefore each term converges:
\begin{equation}\label{eq:gradient-convergence}
 \begin{gathered}
 |\nabla u_n|_2^2\to|\nabla u|_2^2,\qquad
 \int_{\R^N}V_+u_n^2\dx\to\int_{\R^N}V_+u^2\dx,\\
 \lambda_n|u_n|_p^p\to0,\qquad
 \int_{\R^N}(f(u_n)u_n)_-\dx\to\int_{\R^N}(f(u)u)_-\dx.
 \end{gathered}
\end{equation}
With $u_n\to u$ in $L^2(\R^N)$, this gives $\norm{u_n}_V\to\norm{u}_V$.
Since $u_n\weakto u$ in $E_V$, we get $u_n\to u$ in $E_V$.

By \eqref{eq:gradient-convergence} and Scheff\'e's lemma,
$(f(u_n)u_n)_-\to(f(u)u)_-$ in $L^1(\R^N)$. With
\eqref{eq:positive-part-convergence}, $f(u_n)u_n\to f(u)u$ in
$L^1(\R^N)$. By \eqref{eq:F-domination},
\begin{equation*}
 |F(u_n)|\le C\bigl(|u_n|^2+|f(u_n)u_n|\bigr),
\end{equation*}
and the right-hand side converges in $L^1(\R^N)$. So $F(u_n)\to F(u)$ in
$L^1(\R^N)$ by Vitali's theorem. By \eqref{eq:gradient-convergence} and
\eqref{eq:potential-negative-conv}, $J_{\lambda_n}(u_n)\to J_0(u)$. This
proves \eqref{eq:strong-limit}.

Finally, let \(\liminf_nJ_{\lambda_n}(u_n)>0\). Then \(J_0(u)>0=J_0(0)\),
so \(u\ne0\).
\end{proof}

\begin{proof}[Proof of Theorem \ref{thm:existence}]
If \(0<\lambda\le\mu\le1\), then \(J_\lambda\le J_\mu\) on \(X\). So
\(c_\lambda\le c_\mu\). By Lemma~\ref{lem:MP-bounds},
\(m_1\le c_\lambda\le m_2\) for all \(\lambda\in(0,1]\). Therefore
\begin{equation*}
 c_0:=\lim_{\lambda\to0^+}c_\lambda\in[m_1,m_2]
\end{equation*}
exists. Take \(\lambda_n\in(0,1]\) with \(\lambda_n\to0^+\), and let
\(u_n=u_{\lambda_n}\) be given by \eqref{eq:MP-critical}. Then
\(J_{\lambda_n}(u_n)=c_{\lambda_n}\le m_2\). Lemma~\ref{lem:limit} gives,
along a subsequence, a radial finite-energy weak solution \(u_0\) with
\(u_n\to u_0\) in \(E_V\) and
\begin{equation*}
 J_0(u_0)=\lim_{n\to\infty}J_{\lambda_n}(u_n)
       =\lim_{n\to\infty}c_{\lambda_n}=c_0\ge m_1>0.
\end{equation*}
In particular, \(u_0\ne0\).
\end{proof}

\section{Proof of Theorem \ref{thm:multiplicity}}\label{sec:multiplicity}

We use a symmetric mountain pass argument. The argument is standard. We
give the details, since all bounds must be uniform in $\lambda$.
Throughout this section we assume \eqref{V1}, and $X$ is given by
\eqref{eq:X} with $1<p<2$ fixed.

We first fix a basis. The space $X$ is separable, and
$C_{c,\rad}^\infty(\R^N)$ is dense in $X$. So we can choose a linearly
independent sequence in $C_{c,\rad}^\infty(\R^N)$ whose span is dense in
$X$. We apply Gram--Schmidt in $E_V$. This gives
$e_1,e_2,\ldots\in C_{c,\rad}^\infty(\R^N)$, orthonormal for the inner
product of $E_V$, with the same span. This span is dense in $X$. Since
$X$ is dense in $E_V$, it is also dense in $E_V$. For $k\ge1$, set
\begin{equation}\label{eq:YZ}
 \begin{split}
 Y_k&=\spanof\{e_1,\ldots,e_k\},\\
 Z_k&=\left\{u\in X:
       \int_{\R^N}(\nabla u\cdot\nabla e_i+(1+V_+(x))ue_i)\dx=0,
                                  \ 1\le i\le k\right\},
 \end{split}
\end{equation}
and $Z_0=X$. We also recall the equivariant deformation lemma; see
\cite{Willem}.

\begin{lemma}\label{lem:deformation}
Let \(I\in C^1(X,\R)\) be even and satisfy the Palais--Smale
condition. Suppose that \(c>0\) is not a critical value of \(I\).
Then there exist \(0<\varepsilon<c/2\) and an odd continuous map
\(\eta:X\to X\) such that
\begin{enumerate}
\item[(1)] $\eta(u)=u$ if $I(u)\notin[c-2\varepsilon,c+2\varepsilon]$;
\item[(2)] $I(\eta(u))\le I(u)$ for all $u\in X$;
\item[(3)] $\eta(\{I\le c+\varepsilon\})\subset\{I\le c-\varepsilon\}$.
\end{enumerate}
\end{lemma}

The next lemma is a version of the intersection lemma used in the proof
of the fountain theorem; compare \cite{Willem}. We include the short
proof.

\begin{lemma}\label{lem:intersection}
Let $k\ge1$ and $\rho>0$. Let $h:\overline B^k\to X$ be continuous and
odd, where $\overline B^k$ is the closed unit ball of $\R^k$. If
$\norm{h(y)}_V>\rho$ for every $y\in S^{k-1}$, then there is
$y\in\overline B^k$ with
\begin{equation}\label{eq:intersection}
 h(y)\in Z_{k-1},\qquad \norm{h(y)}_V=\rho.
\end{equation}
\end{lemma}

\begin{proof}
Let $\Omega=\{y\in B^k:\norm{h(y)}_V<\rho\}$, where $B^k$ is the open
unit ball of $\R^k$. Since $h$ is odd and continuous, $h(0)=0$, and
$\Omega$ is a symmetric open neighborhood of $0$. As
$\norm{h(y)}_V>\rho$ on $S^{k-1}$, we have $\overline\Omega\subset B^k$.
Hence $\norm{h(y)}_V=\rho$ on $\partial\Omega$. If $k=1$, then $Z_0=X$,
and any point of $\partial\Omega$ gives \eqref{eq:intersection}.

Let $k\ge2$, and let $\Lambda:X\to\R^{k-1}$ map $u$ to the vector of the
first $k-1$ inner products in \eqref{eq:YZ}. Then $\Lambda$ is linear
and continuous, and $\ker\Lambda=Z_{k-1}$. Suppose that
\eqref{eq:intersection} fails. Then $\Lambda\circ h$ maps
$\partial\Omega$ into $\R^{k-1}\setminus\{0\}$, and it is odd and
continuous. This contradicts the Borsuk--Ulam theorem: every odd
continuous map from the boundary of a bounded symmetric open
neighborhood of $0$ in $\R^k$ into $\R^{k-1}$ has a zero.
\end{proof}

\begin{lemma}\label{lem:tau}
Fix $q\in(2,2^*)$ and let
\begin{equation*}
 \tau_k=\sup_{\substack{u\in Z_k\\\norm{u}_V=1}}|u|_q,
       \qquad k\ge0.
\end{equation*}
Then $0<\tau_k<\infty$ and $\tau_k\to0$ as $k\to\infty$.
\end{lemma}

\begin{proof}
By the Sobolev inequality, $\tau_k<\infty$. Since $Z_k$ has finite
codimension in $X$, $Z_k\ne\{0\}$, and so $\tau_k>0$. Since
$Z_{k+1}\subset Z_k$, $\tau_k$ is nonincreasing. Suppose that
$\tau_k\ge\delta>0$ for all $k$. Then there are $w_k\in Z_k$ with
$\norm{w_k}_V=1$ and $|w_k|_q\ge\delta/2$. We may assume $w_k\weakto w$
in $E_V$. For each fixed $i$, the inner product of $w_k$ with $e_i$ in
$E_V$ vanishes for $k\ge i$. So $w$ is orthogonal to every $e_i$. The
span of the $e_i$ is dense in $E_V$, and so $w=0$. The embedding
$E_V\hookrightarrow L^q(\R^N)$ is compact, since
$\norm{u}\le\norm{u}_V$ and $H^1_{\rad}(\R^N)\hookrightarrow L^q(\R^N)$
is compact. Hence $|w_k|_q\to0$, a contradiction.
\end{proof}

\begin{lemma}\label{lem:levels}
There exist an integer \(j_0\ge1\) and numbers
\(0<a_j\le b_j<\infty\) \((j\ge j_0)\), independent of
\(\lambda\in(0,1]\), with \(a_j\to\infty\) as \(j\to\infty\), such that
for every \(j\ge j_0\) and every \(\lambda\in(0,1]\) there is
\(u_{\lambda,j}\in X\) satisfying
\begin{equation*}
 J_\lambda'(u_{\lambda,j})=0,\qquad
       a_j\le J_\lambda(u_{\lambda,j})\le b_j.
\end{equation*}
\end{lemma}

\begin{proof}
\emph{Step 1: the energy is negative far out in $Y_j$.} For every fixed
\(j\),
\begin{equation}\label{eq:finite-negative}
 \sup_{\substack{u\in Y_j,\ \norm{u}_V=t\\0<\lambda\le1}}
       J_\lambda(u)\longrightarrow-\infty
       \quad\text{as}\quad t\to\infty.
\end{equation}
Indeed, the unit sphere of \(Y_j\) is compact. On it, \(|w|_2\) has a
positive minimum and \(|w|_p\) has a finite maximum. Let \(w\in Y_j\)
with \(\norm{w}_V=1\). By \eqref{eq:F-lower} and $\lambda\le1$,
\begin{equation*}
 J_\lambda(tw)
 \le t^2\left(\frac12-L|w|_2^2\right)
       +\left(C_L+\frac1p\right)t^p|w|_p^p.
\end{equation*}
For \(L\) large, the coefficient of \(t^2\) is at most \(-1\) on this
sphere. Since \(p<2\), \eqref{eq:finite-negative} follows.

\emph{Step 2: the energy is large on spheres in $Z_{j-1}$.} Let
$\tau_k$ be given by Lemma~\ref{lem:tau}, and put
\(r_j=\tau_{j-1}^{-1}\). Then \(|u|_q\le1\) for
\(u\in Z_{j-1}\) with \(\norm{u}_V=r_j\). By \eqref{eq:lower-geometry},
\begin{equation*}
 \inf_{\substack{u\in Z_{j-1}\\\norm{u}_V=r_j}}J_\lambda(u)
       \ge\frac14r_j^2-C_0=:a_j\longrightarrow+\infty,
\end{equation*}
uniformly in \(\lambda\). Fix \(j_0\ge1\) with \(a_j>0\) for all
\(j\ge j_0\).

\emph{Step 3: minimax values.} Let \(j\ge j_0\). By
\eqref{eq:finite-negative}, there is \(\rho_j>r_j\) such that
\(J_\lambda(u)<0\) for all \(u\in Y_j\) with \(\norm{u}_V=\rho_j\) and
all \(0<\lambda\le1\). Define
\begin{equation*}
 \begin{split}
 \Gamma_j&=\Bigl\{h\in C(\overline B^j,X):h\text{ is odd and }
             h(y)=\rho_j\sum_{i=1}^jy_ie_i\text{ for }y\in S^{j-1}\Bigr\},\\
 c_{\lambda,j}&=\inf_{h\in\Gamma_j}\max_{y\in\overline B^j}J_\lambda(h(y)).
 \end{split}
\end{equation*}
Let \(h\in\Gamma_j\). Since the \(e_i\) are orthonormal in \(E_V\),
\(\norm{h(y)}_V=\rho_j>r_j\) for \(y\in S^{j-1}\). By
Lemma~\ref{lem:intersection} with \(k=j\) and \(\rho=r_j\), there is
\(y\) with \(h(y)\in Z_{j-1}\) and \(\norm{h(y)}_V=r_j\). Hence
\(c_{\lambda,j}\ge a_j\). The map
\(h_0(y)=\rho_j\sum_{i=1}^jy_ie_i\) belongs to \(\Gamma_j\), and
\(J_\lambda\le J_1\). So
\begin{equation*}
 c_{\lambda,j}\le\max_{y\in\overline B^j}J_1(h_0(y))=:b_j<\infty,
\end{equation*}
since \(h_0(\overline B^j)\) is compact in \(X\).

\emph{Step 4: $c_{\lambda,j}$ is a critical value.} Fix
\(\lambda\in(0,1]\). Suppose that \(c_{\lambda,j}\) is not a critical
value of \(J_\lambda\). Lemmas~\ref{lem:PS} and~\ref{lem:deformation}
give \(\varepsilon\in(0,c_{\lambda,j}/2)\) and an odd map \(\eta\).
Choose \(h\in\Gamma_j\) with
\begin{equation*}
 \max_{y\in\overline B^j}J_\lambda(h(y))\le c_{\lambda,j}+\varepsilon.
\end{equation*}
For \(y\in S^{j-1}\), \(J_\lambda(h(y))<0<c_{\lambda,j}-2\varepsilon\).
So \(\eta(h(y))=h(y)\). Thus \(\eta\circ h\in\Gamma_j\), but
\begin{equation*}
 \max_{y\in\overline B^j}J_\lambda(\eta(h(y)))
       \le c_{\lambda,j}-\varepsilon.
\end{equation*}
This contradicts the definition of \(c_{\lambda,j}\). So there is
\(u_{\lambda,j}\in X\) with \(J_\lambda'(u_{\lambda,j})=0\) and
\(J_\lambda(u_{\lambda,j})=c_{\lambda,j}\in[a_j,b_j]\).
\end{proof}

\begin{proof}[Proof of Theorem \ref{thm:multiplicity}]
Let \(u_{\lambda,j}\), \(a_j\) and \(b_j\) be as in
Lemma~\ref{lem:levels}. Since \(a_j\to\infty\) and every \(b_j\) is
finite, we can choose \(j_0\le j_1<j_2<\cdots\) such that
\begin{equation}\label{eq:separated-levels}
 a_{j_{m+1}}>b_{j_m}\qquad(m\ge1).
\end{equation}
Fix \(\lambda_n\in(0,1]\) with \(\lambda_n\to0^+\). For each \(m\), we
apply Lemma~\ref{lem:limit} to \(\{u_{\lambda_n,j_m}\}_n\). Along a
subsequence, which may depend on \(m\), we obtain a radial finite-energy
weak solution \(u_m\) with
\begin{equation}\label{eq:level-survival}
 u_{\lambda_n,j_m}\to u_m\text{ in }E_V,\qquad
       a_{j_m}\le J_0(u_m)
           =\lim_{n\to\infty}J_{\lambda_n}(u_{\lambda_n,j_m})
           \le b_{j_m}.
\end{equation}
By \eqref{eq:separated-levels} and \eqref{eq:level-survival},
\begin{equation*}
 0<J_0(u_1)<J_0(u_2)<\cdots,\qquad J_0(u_m)\ge a_{j_m}\to\infty.
\end{equation*}
In particular, every \(u_m\) is nonzero, and the pairs \(\{\pm u_m\}\)
are distinct.
\end{proof}

\section{General nonlinearities and the Pohozaev identity}\label{sec:general}

Throughout this section we assume \eqref{V2} and
\eqref{g1}--\eqref{g3}, and fix $1<p<r$. Since $V$ is bounded,
\eqref{V1} holds with $M_0=|V_-|_\infty$. Also $E_V=H^1_{\rad}(\R^N)$,
and
\begin{equation*}
 \norm{u}^2\le\norm{u}_V^2\le(1+|V|_\infty)\norm{u}^2.
\end{equation*}
So the space $X$ in \eqref{eq:X} is $X=H^1_{\rad}(\R^N)\cap L^p(\R^N)$,
with a norm equivalent to $\norm{u}+|u|_p$. As in
Section~\ref{sec:existence}, $X$ is separable and reflexive, and
$C_{c,\rad}^\infty(\R^N)$ is dense in $X$. Since $p<r<2$, the vectors
$e_i$, the spaces $Z_k$ in \eqref{eq:YZ} and Lemmas~\ref{lem:intersection}
and~\ref{lem:tau} are also available here. For $0\le\lambda\le1$, the
functional $I_\lambda$ is given by \eqref{eq:I-intro}. The Pohozaev
functional is
\begin{equation*}
 \begin{split}
 P_\lambda(u)={}&\frac{N-2}{2}\int_{\R^N}|\nabla u|^2\dx
       +\frac12\int_{\R^N}(NV(x)+W(x))u^2\dx\\
 &+\frac{N\lambda}{p}\int_{\R^N}|u|^p\dx
       -N\int_{\R^N}G(u)\dx.
 \end{split}
\end{equation*}

\begin{lemma}\label{lem:g-estimates}
The functionals $I_\lambda$ belong to $C^1(X,\R)$ and are even.
There is $C>0$ such that
\begin{equation}\label{eq:g-growth}
 \begin{split}
 |g(t)|&\le C\bigl(|t|^{r-1}+|t|^{2^*-1}\bigr),\\
 |G(t)|&\le C\bigl(|t|^r+|t|^{2^*}\bigr),\qquad
 G_+(t)\le C|t|^{2^*}.
 \end{split}
\end{equation}
Moreover,
\begin{equation}\label{eq:G-mass}
 \lim_{t\to0}\frac{-G(t)}{t^2}=+\infty.
\end{equation}
For every $A>0$, $\varepsilon>0$ and $q\in(2,2^*)$, there is
$C_{A,\varepsilon,q}>0$ such that
\begin{equation}\label{eq:G-upper}
 G(t)\le -At^2+C_{A,\varepsilon,q}|t|^q
                      +\varepsilon|t|^{2^*}\qquad(t\in\R).
\end{equation}
If $u_n\weakto u$ in $H^1_{\rad}(\R^N)$, then
\begin{equation}\label{eq:Gplus-compact}
 G_+(u_n)\longrightarrow G_+(u)\quad\text{in }L^1(\R^N).
\end{equation}
\end{lemma}

\begin{proof}
The bounds \eqref{eq:g-growth} follow from \eqref{g2}, \eqref{g3} and
the continuity of $g$. Let $A>0$. By \eqref{g1}, for $t>0$ small,
\[
 g(t)\le-2At,\qquad G(t)\le-At^2.
\]
Since $G$ is even, $\liminf_{t\to0}(-G(t)/t^2)\ge A$ for every $A>0$.
This is \eqref{eq:G-mass}. Also $G_+(t)=0$ for $|t|$ small.

Let $A,\varepsilon>0$ and $2<q<2^*$. By \eqref{eq:G-mass},
$\bigl(G(t)+At^2-\varepsilon|t|^{2^*}\bigr)_+=0$ for $|t|$ small. By
\eqref{g3}, for $|t|$ large,
\[
 G(t)\le G_+(t)\le\frac\varepsilon2|t|^{2^*},\qquad
 At^2\le\frac\varepsilon2|t|^{2^*}.
\]
So this function also vanishes for $|t|$ large. By continuity,
\[
 \sup_{t\ne0}
 \frac{\bigl(G(t)+At^2-\varepsilon|t|^{2^*}\bigr)_+}{|t|^q}<\infty,
\]
which is \eqref{eq:G-upper}.

Let $\chi\in C_c(\R)$ satisfy $0\le\chi\le1$, $\chi=1$ on $[-1,1]$ and
$\chi=0$ outside $[-2,2]$. Then
\[
 |\chi(t)g(t)|\le C|t|^{r-1},\qquad
 |(1-\chi(t))g(t)|\le C|t|^{2^*-1}.
\]
Since $p<r<2<2^*$, $X\hookrightarrow L^r(\R^N)\cap L^{2^*}(\R^N)$. Let
$u_n\to u$ in $X$. We apply Vitali's theorem along subsequences that
converge almost everywhere. This gives
\[
 \begin{aligned}
 \chi(u_n)g(u_n)&\to\chi(u)g(u)
       &&\text{in }L^{r/(r-1)}(\R^N),\\
 (1-\chi(u_n))g(u_n)&\to(1-\chi(u))g(u)
       &&\text{in }L^{2^*/(2^*-1)}(\R^N).
 \end{aligned}
\]
By H\"older's inequality, for $h\to0$ in $X$,
\[
 \left|\int_{\R^N}\bigl(G(u+h)-G(u)-g(u)h\bigr)\dx\right|
 \le\int_0^1\left|\int_{\R^N}
       \bigl(g(u+sh)-g(u)\bigr)h\dx\right|\,\mathrm{d}s
       =o\bigl(\norm{h}_X\bigr).
\]
Hence $I_\lambda\in C^1(X,\R)$ and
\[
 \pair{I_\lambda'(u)}{v}
 =\int_{\R^N}(\nabla u\cdot\nabla v+V(x)uv)\dx
 +\lambda\int_{\R^N}|u|^{p-2}uv\dx
       -\int_{\R^N}g(u)v\dx.
\]
Since $g$ is odd, $I_\lambda$ is even.

Now let $u_n\weakto u$ in $H^1_{\rad}(\R^N)$. Then $u_n\to u$ in
$L^q(\R^N)$ for $2<q<2^*$, and $\sup_n|u_n|_{2^*}<\infty$. Let
$\varepsilon>0$. Choose $T>1$ with $G_+(t)\le\varepsilon|t|^{2^*}$ for
$|t|\ge T$. Let $\chi_T\in C_c(\R)$ satisfy $0\le\chi_T\le1$,
$\chi_T=1$ on $[-T,T]$ and $\chi_T=0$ outside $[-2T,2T]$. Then
\[
 0\le(1-\chi_T(t))G_+(t)\le\varepsilon|t|^{2^*},\qquad
 0\le\chi_T(t)G_+(t)\le C_T|t|^q.
\]
By Vitali's theorem, applied along subsequences that converge almost
everywhere, $\chi_T(u_n)G_+(u_n)\to\chi_T(u)G_+(u)$ in $L^1(\R^N)$. It
follows that
\[
 \limsup_{n\to\infty}|G_+(u_n)-G_+(u)|_1
 \le\varepsilon\Bigl(\sup_n|u_n|_{2^*}^{2^*}+|u|_{2^*}^{2^*}\Bigr).
\]
Letting $\varepsilon\to0$, we obtain \eqref{eq:Gplus-compact}.
\end{proof}

\begin{lemma}\label{lem:g-tail}
Let $u_n\in X$, $0\le\lambda_n\le1$, and suppose that
\begin{equation*}
 \sup_n\norm{u_n}<\infty,\qquad
       \sup_n I_{\lambda_n}(u_n)<\infty.
\end{equation*}
Then $\{u_n\}$ is relatively compact in $L^2(\R^N)$. Moreover, let
\begin{equation*}
 m(\rho)=\inf_{0<|t|\le\rho}\frac{-G(t)}{t^2}.
\end{equation*}
Then, for all sufficiently large $R$,
\begin{equation}\label{eq:general-tail}
 \sup_n\int_{|x|>R}|u_n|^2\dx
       \le\frac{C}{m(CR^{-(N-1)/2})}\longrightarrow0.
\end{equation}
\end{lemma}

\begin{proof}
The proof follows that of Lemma~\ref{lem:compactness}. By
\eqref{eq:g-growth} and the Sobolev inequality,
\begin{equation*}
 \begin{split}
 \int_{\R^N}G_-(u_n)\dx+\frac{\lambda_n}{p}|u_n|_p^p
 &=I_{\lambda_n}(u_n)
   -\frac12\int_{\R^N}(|\nabla u_n|^2+V(x)u_n^2)\dx
   +\int_{\R^N}G_+(u_n)\dx\\
 &\le C+\frac{M_0}{2}|u_n|_2^2
                       +C|u_n|_{2^*}^{2^*}\le C.
 \end{split}
\end{equation*}
By \eqref{g1}, $g<0$ on $(0,t_0)$ for some $t_0>0$. So $G_-=-G$ on
$[-t_0,t_0]$. Also $m(\rho)\to\infty$ as $\rho\to0^+$, by
\eqref{eq:G-mass}. By \eqref{eq:radial}, for $R$ large,
$|u_n(x)|\le CR^{-(N-1)/2}\le t_0$ for $|x|>R$. Thus
\[
 G_-(u_n)\ge m(CR^{-(N-1)/2})|u_n|^2\qquad\text{on }\{|x|>R\}.
\]
Integrating over $\{|x|>R\}$, we obtain \eqref{eq:general-tail}. The
weak limit satisfies the same tail bound, by Fatou's lemma. With local
$L^2$ compactness, this gives strong convergence in $L^2(\R^N)$ along a
subsequence.
\end{proof}

\begin{lemma}\label{lem:Pohozaev}
Let $\lambda\ge0$. Suppose that $u\in H^1_{\rad}(\R^N)$,
$G(u)\in L^1(\R^N)$, $u\in L^p(\R^N)$ if $\lambda>0$, and $u$ is a weak
solution of
\begin{equation*}
 -\Delta u+V(x)u+\lambda|u|^{p-2}u=g(u),
\end{equation*}
that is, the corresponding integral identity holds for all
$\varphi\in C_{c,\rad}^\infty(\R^N)$. Then $g(u)u\in L^1(\R^N)$,
\begin{equation}\label{eq:g-Nehari}
 \int_{\R^N}(|\nabla u|^2+V(x)u^2)\dx
       +\lambda\int_{\R^N}|u|^p\dx
       =\int_{\R^N}g(u)u\dx,
\end{equation}
and
\begin{equation}\label{eq:P-zero}
 P_\lambda(u)=0.
\end{equation}
\end{lemma}

\begin{proof}
First we prove \eqref{eq:g-Nehari}. Since $g(t)t\le0$ near zero and $g$
has at most critical growth, $(g(t)t)_+\le C|t|^{2^*}$. So
$(g(u)u)_+\in L^1(\R^N)$. On every ball,
$g(u)\in L^{2^*/(2^*-1)}$. Hence the equation extends to radial test
functions in $H^1(\R^N)$ with compact support. Let $\eta_R$ be the
cut-off functions from Step~3 of the proof of Lemma~\ref{lem:limit}.
Testing with $\eta_Ru$ gives
\begin{equation*}
 \begin{split}
 \int_{\R^N}\eta_R(g(u)u)_-\dx
 ={}&\int_{\R^N}\eta_R(g(u)u)_+\dx
       -\int_{\R^N}\eta_R(|\nabla u|^2+V(x)u^2)\dx\\
 &-\lambda\int_{\R^N}\eta_R|u|^p\dx
       -\int_{\R^N}u\nabla u\cdot\nabla\eta_R\dx.
 \end{split}
\end{equation*}
Since $V$ is bounded and the last integral tends to zero, the right-hand
side is bounded above independently of $R$. By Fatou's lemma,
$(g(u)u)_-\in L^1(\R^N)$. Dominated convergence in the tested equation
gives \eqref{eq:g-Nehari}.

For the Pohozaev identity we use the radial equation on annuli. Write
$u(x)=U(|x|)$. Since $U\in H^1_{\mathrm{loc}}((0,\infty))$, $U$ is
continuous on $(0,\infty)$. Testing with $\varphi(x)=\psi(|x|)$,
$\psi\in C_c^\infty((0,\infty))$, we get
\[
 (s^{N-1}U')'=s^{N-1}\bigl(V(s)U+\lambda|U|^{p-2}U-g(U)\bigr)
 \quad\text{in }\mathcal D'((0,\infty)).
\]
The right-hand side is continuous. Hence $s^{N-1}U'\in C^1((0,\infty))$,
$U\in C^2((0,\infty))$, and
\begin{equation}\label{eq:radial-ODE}
 -U''-\frac{N-1}{s}U'+V(s)U+\lambda|U|^{p-2}U=g(U)\qquad(s>0).
\end{equation}
We multiply \eqref{eq:radial-ODE} by $s^NU'$ and integrate over
$[a,b]\subset(0,\infty)$. We use
\[
 \begin{aligned}
 s^NU'U''&=\Bigl(\frac{s^N}2|U'|^2\Bigr)'-\frac N2s^{N-1}|U'|^2,\\
 s^NV(s)UU'&=\Bigl(\frac{s^N}2V(s)U^2\Bigr)'
       -\frac12s^{N-1}\bigl(NV(s)+sV'(s)\bigr)U^2,
 \end{aligned}
\]
and similar identities for $|U|^{p-2}UU'$ and $g(U)U'$. This gives
\begin{equation}\label{eq:P-annuli}
 \begin{split}
 &\frac{N-2}{2}\int_a^b s^{N-1}|U'|^2\,\mathrm{d}s
 +\frac12\int_a^b s^{N-1}(NV(s)+sV'(s))U^2\,\mathrm{d}s\\
 &\quad+\frac{N\lambda}{p}\int_a^b s^{N-1}|U|^p\,\mathrm{d}s
       -N\int_a^b s^{N-1}G(U)\,\mathrm{d}s\\
 &=\left[s^N\left(-\frac12|U'|^2+\frac12V(s)U^2
                         +\frac\lambda p|U|^p-G(U)\right)\right]_a^b.
 \end{split}
\end{equation}
Note that $sV'(s)=W(x)$ for $|x|=s$. Since $W$ is bounded, all integrals
on the left converge absolutely as $a\to0^+$ and $b\to\infty$. Moreover,
\begin{equation*}
 \int_0^\infty s^{N-1}
   \bigl(|U'|^2+|V(s)|U^2+\lambda|U|^p+|G(U)|\bigr)
       \,\mathrm{d}s<\infty.
\end{equation*}
So there are $a_n\to0^+$ and $b_n\to\infty$ along which $s^N$ times this
integrand tends to zero. Otherwise the integrand would be at least
$cs^{-N}$ for some $c>0$ near $0$ or near $\infty$, and the integral
would dominate $c\int_0^1s^{-1}\,\mathrm{d}s$ or
$c\int_1^\infty s^{-1}\,\mathrm{d}s$. Take $a=a_n$, $b=b_n$ in
\eqref{eq:P-annuli}, multiply by the area of the unit sphere and let
$n\to\infty$. The boundary terms vanish, and \eqref{eq:P-zero} follows.
No regularity of $U$ at the origin is needed.
\end{proof}

\begin{lemma}\label{lem:general-bounds}
For every $M>0$ there is $C_M>0$, independent of $0<\lambda\le1$, such
that
\begin{equation*}
 I_\lambda(u)\le M,\qquad |P_\lambda(u)|\le1,
       \qquad u\in X
\end{equation*}
imply
\begin{equation*}
 \norm{u}^2+\lambda|u|_p^p
       +\int_{\R^N}G_-(u)\dx\le C_M.
\end{equation*}
\end{lemma}

\begin{proof}
We subtract $\frac1NP_\lambda$ from $I_\lambda$.
This removes all terms except the gradient term and the $W$ term. The
smallness of $W_+$ then bounds $|\nabla u|_2$. The infinite-mass
condition turns this into a bound on $|u|_2$.

In detail, we have
\begin{equation}\label{eq:IP-cancellation}
 I_\lambda(u)-\frac1N P_\lambda(u)
 =\frac1N\int_{\R^N}|\nabla u|^2\dx
       -\frac1{2N}\int_{\R^N}W(x)u^2\dx.
\end{equation}
By H\"older's inequality and \eqref{eq:Sobolev},
\begin{equation}\label{eq:W-form-bound}
 \int_{\R^N}W_+(x)u^2\dx
 \le |W_+|_{N/2}|u|_{2^*}^2
 \le\frac{|W_+|_{N/2}}{S}|\nabla u|_2^2.
\end{equation}
Since $W\le W_+$, \eqref{eq:IP-cancellation} and \eqref{eq:W-form-bound}
give
\begin{equation*}
 \frac1N\left(1-\frac{|W_+|_{N/2}}{2S}\right)|\nabla u|_2^2
       \le I_\lambda(u)-\frac1N P_\lambda(u)
       \le M+\frac1N.
\end{equation*}
By \eqref{V2}, $|W_+|_{N/2}<2S$. So $|\nabla u|_2$ and, by
\eqref{eq:Sobolev}, $|u|_{2^*}$ are bounded independently of $\lambda$.

Next we bound $|u|_2$. Since $G_+\ge0$, we have $G_-\ge-G$. So
\eqref{eq:G-mass} gives $\delta>0$ such that
$G_-(t)\ge (M_0+1)t^2$ for $|t|\le\delta$. For $|t|>\delta$,
$t^2\le\delta^{2-2^*}|t|^{2^*}$. Therefore
\begin{equation}\label{eq:mass-absorption}
 |u|_2^2\le \frac{1}{M_0+1}\int_{\R^N}G_-(u)\dx
                         +\delta^{2-2^*}|u|_{2^*}^{2^*}.
\end{equation}
By the definition of $I_\lambda$, \eqref{eq:g-growth} and
\eqref{eq:mass-absorption},
\begin{equation*}
 \begin{split}
 \int_{\R^N}G_-(u)\dx+\frac\lambda p|u|_p^p
 &\le M+\frac{M_0}{2}|u|_2^2+\int_{\R^N}G_+(u)\dx\\
 &\le C_M+\frac{M_0}{2M_0+2}\int_{\R^N}G_-(u)\dx.
 \end{split}
\end{equation*}
Since $M_0/(2M_0+2)<1/2$, $\int_{\R^N}G_-(u)\dx$ and
$\lambda|u|_p^p$ are bounded. By \eqref{eq:mass-absorption}, so is
$|u|_2$.
\end{proof}

\begin{lemma}\label{lem:PSP}
Fix $0<\lambda\le1$. Suppose that $u_n\in X$ satisfies
\begin{equation}\label{eq:PSP}
 I_\lambda(u_n)\to c,\qquad
 I_\lambda'(u_n)\to0\text{ in }X^*,\qquad
 P_\lambda(u_n)\to0.
\end{equation}
Then a subsequence converges in $X$ to a critical point of
$I_\lambda$ at level $c$.
\end{lemma}

\begin{proof}
By Lemma~\ref{lem:general-bounds}, $\{u_n\}$ is bounded in $X$, since
$\lambda>0$ is fixed. By Lemma~\ref{lem:g-tail} and the compactness of
the radial embedding, up to a subsequence,
\[
 u_n\weakto u\ \text{in }X,\qquad
 u_n\to u\ \text{in }L^q(\R^N)\ (2\le q<2^*),\qquad
 u_n\to u\ \text{a.e. in }\R^N.
\]
Let $R>0$. By \eqref{eq:g-growth}, $\{g(u_n)\}$ is bounded in
$L^{2^*/(2^*-1)}(B_R)$. So $g(u_n)\to g(u)$ in $L^1(B_R)$ by Vitali's
theorem. Since $u_n\to u$ in $L^p(B_R)$,
$|u_n|^{p-2}u_n\to|u|^{p-2}u$ in $L^{p/(p-1)}(B_R)$. Therefore, for
$\varphi\in C_{c,\rad}^\infty(\R^N)$,
\[
 \pair{I_\lambda'(u)}{\varphi}
 =\lim_{n\to\infty}\pair{I_\lambda'(u_n)}{\varphi}=0.
\]
By density, $I_\lambda'(u)=0$. Lemma~\ref{lem:Pohozaev} gives
$P_\lambda(u)=0$. As $V,W\in L^\infty(\R^N)$ and $u_n\to u$ in
$L^2(\R^N)$,
\begin{equation}\label{eq:VW-convergence}
 \int_{\R^N}(NV(x)+W(x))u_n^2\dx
       \longrightarrow\int_{\R^N}(NV(x)+W(x))u^2\dx.
\end{equation}
Since $P_\lambda(u_n)\to0$,
\[
 \begin{split}
 \frac{N-2}{2}|\nabla u_n|_2^2+\frac{N\lambda}{p}|u_n|_p^p
 ={}&-\frac12\int_{\R^N}(NV(x)+W(x))u_n^2\dx\\
  &+N\int_{\R^N}\bigl(G_+(u_n)-G_-(u_n)\bigr)\dx+o(1).
 \end{split}
\]
By \eqref{eq:VW-convergence}, \eqref{eq:Gplus-compact}, Fatou's lemma
for $G_-$ and $P_\lambda(u)=0$,
\[
 \begin{aligned}
 &\limsup_{n\to\infty}\left(
       \frac{N-2}{2}|\nabla u_n|_2^2
       +\frac{N\lambda}{p}|u_n|_p^p\right)\\
 &\quad\le-\frac12\int_{\R^N}(NV(x)+W(x))u^2\dx
       +N\int_{\R^N}G_+(u)\dx
       -N\int_{\R^N}G_-(u)\dx\\
 &\quad=\frac{N-2}{2}|\nabla u|_2^2
       +\frac{N\lambda}{p}|u|_p^p
 \le\liminf_{n\to\infty}\left(
       \frac{N-2}{2}|\nabla u_n|_2^2
       +\frac{N\lambda}{p}|u_n|_p^p\right).
 \end{aligned}
\]
Each of the two terms is weakly lower semicontinuous. So
$|\nabla u_n|_2\to|\nabla u|_2$ and $|u_n|_p\to|u|_p$. The spaces
$L^2(\R^N)$ and $L^p(\R^N)$ are uniformly convex. So weak convergence
and convergence of norms give $\nabla u_n\to\nabla u$ in $L^2(\R^N)$
and $u_n\to u$ in $L^p(\R^N)$. With $u_n\to u$ in $L^2(\R^N)$, this
gives $u_n\to u$ in $X$. Finally, $I_\lambda(u)=\lim_nI_\lambda(u_n)=c$.
\end{proof}

\begin{proof}[Proof of Theorem~\ref{thm:general-limit}]
We may assume $\lambda_n\le1$. Since $u_n\in X$, \eqref{eq:g-growth}
gives $G(u_n)\in L^1(\R^N)$. So Lemmas~\ref{lem:Pohozaev}
and~\ref{lem:general-bounds} give $P_{\lambda_n}(u_n)=0$ and
\begin{equation}\label{eq:zero-limit-bound}
 \norm{u_n}^2+\lambda_n|u_n|_p^p
       +\int_{\R^N}G_-(u_n)\dx\le C.
\end{equation}
By Lemma~\ref{lem:g-tail} and the compactness of the radial embedding,
up to a subsequence,
\begin{equation*}
 u_n\weakto u\text{ in }H^1(\R^N),\qquad
 u_n\to u\text{ in }L^q(\R^N)\ (2\le q<2^*),\qquad
 u_n\to u\text{ a.e. in }\R^N.
\end{equation*}

\emph{The limit is a solution.} As in the proof of Lemma~\ref{lem:PSP},
$g(u_n)\to g(u)$ in $L^1_{\mathrm{loc}}(\R^N)$. The functions $u_n$ and
$V$ are radial. So the averaging argument in Step~2 of the proof of
Lemma~\ref{lem:limit} shows that the equation for $u_n$ holds for every
$\varphi\in C_c^\infty(\R^N)$. For such $\varphi$,
\begin{equation*}
 \left|\lambda_n\int_{\R^N}|u_n|^{p-2}u_n\varphi\dx\right|
 \le\lambda_n^{1/p}(\lambda_n|u_n|_p^p)^{(p-1)/p}|\varphi|_p
       \longrightarrow0.
\end{equation*}
Hence $u$ satisfies \eqref{eq:general-weak}. By Fatou's lemma and
\eqref{eq:zero-limit-bound}, $G_-(u)\in L^1(\R^N)$. By
\eqref{eq:g-growth}, $G_+(u)\in L^1(\R^N)$. Thus $u$ is a radial
finite-energy weak solution of \eqref{eq:general}.
Lemma~\ref{lem:Pohozaev} with $\lambda=0$ gives $g(u)u\in L^1(\R^N)$ and
$P_0(u)=0$.

\emph{Strong convergence.} The potential integrals converge as in
\eqref{eq:VW-convergence}. By $P_{\lambda_n}(u_n)=0$, $P_0(u)=0$,
\eqref{eq:Gplus-compact} and Fatou's lemma for $G_-$,
\begin{equation*}
 \limsup_{n\to\infty}
 \left(\frac{N-2}{2}|\nabla u_n|_2^2
                         +\frac{N\lambda_n}{p}|u_n|_p^p\right)
 \le\frac{N-2}{2}|\nabla u|_2^2
 \le\liminf_{n\to\infty}\frac{N-2}{2}|\nabla u_n|_2^2.
\end{equation*}
So $u_n\to u$ in $H^1(\R^N)$ and $\lambda_n|u_n|_p^p\to0$. We rearrange
$P_{\lambda_n}(u_n)=0$ and get
\begin{equation*}
 \begin{split}
 \int_{\R^N}G_-(u_n)\dx
 ={}&\int_{\R^N}G_+(u_n)\dx
       -\frac{N-2}{2N}|\nabla u_n|_2^2\\
 &-\frac1{2N}\int_{\R^N}(NV(x)+W(x))u_n^2\dx
       -\frac{\lambda_n}{p}|u_n|_p^p
       \longrightarrow\int_{\R^N}G_-(u)\dx.
 \end{split}
\end{equation*}
The limit follows from $P_0(u)=0$. Since $G_-(u_n)\ge0$ and
$G_-(u_n)\to G_-(u)$ a.e., Scheff\'e's lemma gives $G_-(u_n)\to G_-(u)$
in $L^1(\R^N)$. With \eqref{eq:Gplus-compact}, $G(u_n)\to G(u)$ in
$L^1(\R^N)$. Since $V$ is bounded, $u_n\to u$ in $H^1(\R^N)$ and
$\lambda_n|u_n|_p^p\to0$, we get $I_{\lambda_n}(u_n)\to I_0(u)$.
Finally, if $\liminf_nI_{\lambda_n}(u_n)>0$, then $I_0(u)>0$, so
$u\ne0$.
\end{proof}

Next we build odd maps on which $I_\lambda$ is negative.

\begin{lemma}\label{lem:negative-maps}
Assume also \eqref{g4}. There is $\lambda_0\in(0,1]$ such that, for
every integer $k\ge1$ and every $R_0>0$, there is an odd continuous map
$\gamma:S^{k-1}\to X$ with
\begin{equation}\label{eq:negative-maps}
 \inf_{\omega\in S^{k-1}}\norm{\gamma(\omega)}_V>R_0,\qquad
 \sup_{\substack{\omega\in S^{k-1}\\0\le\lambda\le\lambda_0}}
             I_\lambda(\gamma(\omega))<0.
\end{equation}
The constant $\lambda_0$ does not depend on $k$ and $R_0$.
\end{lemma}

\begin{proof}
The functions $\gamma(\omega)$ take values close
to $\pm\xi$ on an annulus, and then they are dilated. Under
$x\mapsto x/T$, the gradient term scales like $T^{N-2}$, and all other
terms of $I_\lambda$ scale like $T^N$.

\emph{Step 1.} By \eqref{g4}, we can choose $\lambda_0\in(0,1]$ with
\[
 \frac{\lambda_0}{p}\xi^p<G(\xi)-\frac12\Bigl(\sup_{\R^N}V\Bigr)\xi^2.
\]
Here $\lambda_0$ depends only on $\xi$, $G(\xi)$, $\sup_{\R^N}V$ and
$p$. Set
\[
 \Psi(t)=G(t)-\frac12\Bigl(\sup_{\R^N}V\Bigr)t^2-\frac{\lambda_0}{p}|t|^p .
\]
Then $\Psi$ is even, $\Psi(0)=0$ and $\Psi(\xi)>0$. Since
$V\le\sup_{\R^N}V$,
\begin{equation}\label{eq:negmap-upper}
 I_\lambda(u)\le\frac12\int_{\R^N}|\nabla u|^2\dx
       -\int_{\R^N}\Psi(u)\dx,
 \qquad u\in X,\ 0\le\lambda\le\lambda_0 .
\end{equation}

\emph{Step 2.} Fix $k\ge1$, and let
$\mathcal A=\{x\in\R^N:1<|x|<2\}$. For $n\ge3$, choose a radial
$\zeta_n\in C_c^\infty(\mathcal A)$ with $0\le\zeta_n\le1$ and
$\zeta_n=1$ on $\{1+1/n\le|x|\le2-1/n\}$. For
$\omega=(\omega_1,\ldots,\omega_k)\in S^{k-1}$, let
\[
 \Pi_\omega(s)=\sum_{i=1}^k\omega_is^{i-1},\qquad
 w_{\omega,n}(x)=\xi\zeta_n(x)\tanh\bigl(n\Pi_\omega(|x|)\bigr).
\]
Then $w_{\omega,n}\in C_{c,\rad}^\infty(\R^N)$, $|w_{\omega,n}|\le\xi$,
$w_{\omega,n}=0$ outside $\mathcal A$, and $w_{-\omega,n}=-w_{\omega,n}$.
Since $w_{\omega,n}$ and $\nabla w_{\omega,n}$ are continuous in
$(\omega,x)\in S^{k-1}\times\overline{\mathcal A}$, the map
$\omega\mapsto w_{\omega,n}$ is continuous from $S^{k-1}$ into $X$.

We claim that there is $n\ge3$ such that
\begin{equation}\label{eq:negmap-positive}
 \int_{\R^N}\Psi(w_{\omega,n})\dx\ge\frac12\Psi(\xi)|\mathcal A|
 \qquad\text{for all }\omega\in S^{k-1}.
\end{equation}
If not, then for every $n\ge3$ there is $\omega_n\in S^{k-1}$ for which
\eqref{eq:negmap-positive} fails. We may assume
$\omega_n\to\omega\in S^{k-1}$. Since $\Pi_\omega$ is a nonzero
polynomial, $\{x\in\mathcal A:\Pi_\omega(|x|)=0\}$ is a finite union of
spheres. For every other $x\in\mathcal A$, $\zeta_n(x)\to1$ and
$n\Pi_{\omega_n}(|x|)\to\pm\infty$. So $w_{\omega_n,n}(x)\to\pm\xi$.
Since $\Psi$ is even, $\Psi(w_{\omega_n,n}(x))\to\Psi(\xi)$. Moreover,
$|\Psi(w_{\omega_n,n})|\le\max_{|t|\le\xi}|\Psi(t)|$ on $\mathcal A$ and
$\Psi(w_{\omega_n,n})=0$ outside $\mathcal A$. By dominated convergence,
\[
 \int_{\R^N}\Psi(w_{\omega_n,n})\dx\longrightarrow\Psi(\xi)|\mathcal A|,
\]
a contradiction. Fix such an $n$ and write $w_\omega=w_{\omega,n}$.

\emph{Step 3.} By \eqref{eq:negmap-positive} and $\Psi(0)=0$,
$w_\omega\ne0$ for every $\omega$. By continuity and compactness,
$\max_{\omega}|\nabla w_\omega|_2<\infty$ and
$\min_{\omega}|w_\omega|_2>0$. Since $N>N-2$, we can choose $T>0$ with
\begin{equation}\label{eq:negmap-T}
 T^{N-2}\max_{\omega\in S^{k-1}}|\nabla w_\omega|_2^2
       <T^N\Psi(\xi)|\mathcal A|,
 \qquad
 T^{N/2}\min_{\omega\in S^{k-1}}|w_\omega|_2>R_0 .
\end{equation}
Define $\gamma(\omega)(x)=w_\omega(x/T)$. The map $u\mapsto u(\cdot/T)$
is linear and bounded on $X$. So $\gamma:S^{k-1}\to X$ is odd and
continuous. Changing variables $x=Ty$, we get
\[
 \begin{gathered}
 |\nabla\gamma(\omega)|_2^2=T^{N-2}|\nabla w_\omega|_2^2,\qquad
 |\gamma(\omega)|_2=T^{N/2}|w_\omega|_2,\\
 \int_{\R^N}\Psi(\gamma(\omega))\dx=T^N\int_{\R^N}\Psi(w_\omega)\dx.
 \end{gathered}
\]
Hence, by \eqref{eq:negmap-upper}, \eqref{eq:negmap-positive} and
\eqref{eq:negmap-T}, for $\omega\in S^{k-1}$ and
$0\le\lambda\le\lambda_0$,
\[
 I_\lambda(\gamma(\omega))
 \le\frac{T^{N-2}}2\max_{\omega'\in S^{k-1}}|\nabla w_{\omega'}|_2^2
       -\frac{T^N}{2}\Psi(\xi)|\mathcal A|<0,
 \qquad
 \norm{\gamma(\omega)}_V\ge|\gamma(\omega)|_2>R_0 .
\]
This proves \eqref{eq:negative-maps}.
\end{proof}

The next lemma adapts the augmented functional method of \cite{HIT} to
$I_\lambda$. We lift $I_\lambda$ to $\R\times X$ by dilations. The
$\theta$-derivative of the lifted functional is $P_\lambda$. So a
deformation argument for the lifted functional gives Palais--Smale
sequences with $P_\lambda\to0$.

\begin{lemma}\label{lem:augmented-minimax}
Let $0<\lambda\le\lambda_0$, and let $\gamma:S^{k-1}\to X$ be an odd
continuous map given by Lemma~\ref{lem:negative-maps}. Set
\begin{equation*}
 \begin{split}
 \Gamma&=\{h\in C(\overline B^k,X):
           h\text{ is odd},\ h=\gamma\text{ on }S^{k-1}\},\\
 c&=\inf_{h\in\Gamma}\max_{y\in\overline B^k}I_\lambda(h(y)).
 \end{split}
\end{equation*}
If $\Gamma\ne\varnothing$ and $c>0$, then there is a sequence satisfying
\eqref{eq:PSP} at level $c$. In particular, $c$ is a critical value of
$I_\lambda$.
\end{lemma}

\begin{proof}
For $\theta\in\R$, let $T_\theta u(x)=u(\mathrm e^{-\theta}x)$. Then
\begin{equation*}
 \begin{gathered}
 |T_\theta u|_2^2=\mathrm e^{N\theta}|u|_2^2,\qquad
 |\nabla T_\theta u|_2^2=\mathrm e^{(N-2)\theta}|\nabla u|_2^2,\\
 |T_\theta u|_p^p=\mathrm e^{N\theta}|u|_p^p,\qquad
 \sup_{|\theta|\le1}\bigl(\norm{T_\theta}_{\mathcal L(X)}
              +\norm{T_{-\theta}}_{\mathcal L(X)}\bigr)<\infty.
 \end{gathered}
\end{equation*}
By density of $C_{c,\rad}^\infty(\R^N)$, $T_\theta u\to u$ in $X$ as
$\theta\to0$.

Let $\mathcal E=\R\times X$ with
$\norm{(\theta,u)}_{\mathcal E}=|\theta|+\norm{u}_X$, and
\begin{equation*}
 \begin{split}
 \widetilde I_\lambda(\theta,u)
 &=I_\lambda(T_\theta u)\\
 &=\frac{\mathrm e^{(N-2)\theta}}2
                      \int_{\R^N}|\nabla u|^2\dx
   +\frac{\mathrm e^{N\theta}}2
                      \int_{\R^N}V(\mathrm e^\theta x)u^2\dx\\
 &\quad+\frac{\lambda\mathrm e^{N\theta}}p
                      \int_{\R^N}|u|^p\dx
   -\mathrm e^{N\theta}\int_{\R^N}G(u)\dx.
 \end{split}
\end{equation*}
Since $\frac{\mathrm d}{\mathrm d\theta}V(\mathrm e^\theta x)
=W(\mathrm e^\theta x)$ and $W$ is bounded, we get
\begin{equation*}
 \begin{split}
 \partial_\theta\widetilde I_\lambda(\theta,u)
 ={}&\frac{N-2}{2}\mathrm e^{(N-2)\theta}
                    \int_{\R^N}|\nabla u|^2\dx
 +\frac{\mathrm e^{N\theta}}2\int_{\R^N}
            \bigl(NV(\mathrm e^\theta x)+W(\mathrm e^\theta x)\bigr)u^2\dx\\
 &+\frac{N\lambda}{p}\mathrm e^{N\theta}
                    \int_{\R^N}|u|^p\dx
   -N\mathrm e^{N\theta}\int_{\R^N}G(u)\dx,
 \end{split}
\end{equation*}
and
\begin{equation*}
 \begin{split}
 \pair{\partial_u\widetilde I_\lambda(\theta,u)}{v}
 ={}&\mathrm e^{(N-2)\theta}\int_{\R^N}\nabla u\cdot\nabla v\dx
   +\mathrm e^{N\theta}\int_{\R^N}V(\mathrm e^\theta x)uv\dx\\
 &+\lambda\mathrm e^{N\theta}\int_{\R^N}|u|^{p-2}uv\dx
   -\mathrm e^{N\theta}\int_{\R^N}g(u)v\dx.
 \end{split}
\end{equation*}
These partial derivatives are continuous in $(\theta,u)$ jointly. For
the terms with $V$ and $W$ we use dominated convergence; for example,
\begin{equation*}
 \left|\int_{\R^N}\bigl(V(\mathrm e^\theta x)-V(\mathrm e^{\theta'}x)\bigr)uv\dx\right|
 \le\left(\int_{\R^N}\bigl|V(\mathrm e^\theta x)-V(\mathrm e^{\theta'}x)\bigr|^2u^2\dx\right)^{1/2}|v|_2,
\end{equation*}
and the first factor tends to zero as $\theta'\to\theta$, since $V$ is
bounded and continuous. For the other terms we use
Lemma~\ref{lem:g-estimates}. Hence
$\widetilde I_\lambda\in C^1(\mathcal E,\R)$. A change of variables
shows that
\begin{equation}\label{eq:augmented-derivative}
 \partial_\theta\widetilde I_\lambda(\theta,u)
     =P_\lambda(T_\theta u),\qquad
 \pair{\partial_u\widetilde I_\lambda(\theta,u)}{v}
     =\pair{I_\lambda'(T_\theta u)}{T_\theta v}.
\end{equation}
Since $g$ is odd, $\widetilde I_\lambda\circ\sigma=\widetilde I_\lambda$
for the linear isometry $\sigma(\theta,u)=(\theta,-u)$ of $\mathcal E$.

\emph{Step 1: the lifted minimax value.} Let
\begin{equation*}
 \widetilde\Gamma=\bigl\{\widetilde h\in
       C(\overline B^k,\mathcal E):\
       \widetilde h(-y)=\sigma\widetilde h(y),\
 \widetilde h(\omega)=(0,\gamma(\omega))
                   \text{ for }\omega\in S^{k-1}\bigr\}.
\end{equation*}
Let $\widetilde h=(\Theta,\Phi)\in\widetilde\Gamma$. Then
$h(y)=T_{\Theta(y)}\Phi(y)$ belongs to $\Gamma$. Indeed,
\begin{equation*}
 h(-y)=T_{\Theta(y)}(-\Phi(y))=-h(y),\qquad
 h(\omega)=T_0\gamma(\omega)=\gamma(\omega).
\end{equation*}
Also $h$ is continuous, because $(\theta,u)\mapsto T_\theta u$ is
continuous from $\R\times X$ into $X$. This continuity follows from
$T_\theta T_{\theta'}=T_{\theta+\theta'}$, the local uniform bound on
$\norm{T_\theta}_{\mathcal L(X)}$ and the continuity of
$\theta\mapsto T_\theta u$. Conversely, $h\mapsto(0,h)$ maps $\Gamma$
into $\widetilde\Gamma$. As
$\widetilde I_\lambda(\widetilde h(y))=I_\lambda(h(y))$,
\begin{equation}\label{eq:equal-minimax}
 \inf_{\widetilde h\in\widetilde\Gamma}
      \max_{y\in\overline B^k}
      \widetilde I_\lambda(\widetilde h(y))=c.
\end{equation}

\emph{Step 2: almost critical points near almost optimal sets.} Take
$\varepsilon_n\to0^+$ with $\varepsilon_n<\min\{c/4,1/16\}$, and set
$\delta_n=\sqrt{\varepsilon_n}$. Then $2\delta_n<1$. Choose
$h_n\in\Gamma$ with
\begin{equation*}
 \max_{y\in\overline B^k} I_\lambda(h_n(y))\le c+\varepsilon_n,
 \qquad
 A_n=\{(0,h_n(y)):y\in\overline B^k\}.
\end{equation*}
The set $A_n$ is compact and invariant under $\sigma$. We claim that
there is $z_n=(\theta_n,v_n)\in\mathcal E$ with
\begin{equation}\label{eq:shadowing}
 \operatorname{dist}_{\mathcal E}(z_n,A_n)\le2\delta_n,\qquad
 |\widetilde I_\lambda(z_n)-c|\le2\varepsilon_n,
 \qquad
 \norm{\widetilde I_\lambda'(z_n)}_{\mathcal E^*}\le8\delta_n.
\end{equation}
Suppose not. Then
$\norm{\widetilde I_\lambda'(z)}_{\mathcal E^*}>8\delta_n$ on the closed
set
\begin{equation*}
 \mathcal N_n=\bigl\{z\in\mathcal E:
 \operatorname{dist}_{\mathcal E}(z,A_n)\le2\delta_n,\
 |\widetilde I_\lambda(z)-c|\le2\varepsilon_n\bigr\}.
\end{equation*}
That is, $\mathcal N_n$ is contained in the open set
$\mathcal O_n=\{z:\norm{\widetilde I_\lambda'(z)}_{\mathcal E^*}>8\delta_n\}$.
On $\mathcal O_n$, $\widetilde I_\lambda$ has a pseudo-gradient vector
field $\mathcal V$ (see \cite{Willem}). It is locally Lipschitz, and
\begin{equation*}
 \norm{\mathcal V(z)}_{\mathcal E}\le2\norm{\widetilde I_\lambda'(z)}_{\mathcal E^*},\qquad
 \pair{\widetilde I_\lambda'(z)}{\mathcal V(z)}\ge\norm{\widetilde I_\lambda'(z)}_{\mathcal E^*}^2.
\end{equation*}
Let $\mathcal V_n=\mathcal V/\norm{\mathcal V}_{\mathcal E}$. It is
locally Lipschitz on $\mathcal O_n$, and
\begin{equation*}
 \norm{\mathcal V_n(z)}_{\mathcal E}=1,\qquad
 \pair{\widetilde I_\lambda'(z)}{\mathcal V_n(z)}
 \ge\frac12\norm{\widetilde I_\lambda'(z)}_{\mathcal E^*}>4\delta_n.
\end{equation*}
We replace $\mathcal V_n(z)$ by
$\tfrac12\bigl(\mathcal V_n(z)+\sigma\mathcal V_n(\sigma z)\bigr)$. Then
$\mathcal V_n(\sigma z)=\sigma\mathcal V_n(z)$ and
$\norm{\mathcal V_n}_{\mathcal E}\le1$. Since
$\widetilde I_\lambda\circ\sigma=\widetilde I_\lambda$, we have
$\widetilde I_\lambda'(\sigma z)=\widetilde I_\lambda'(z)\circ\sigma$.
So
$\pair{\widetilde I_\lambda'(z)}{\sigma\mathcal V_n(\sigma z)}
=\pair{\widetilde I_\lambda'(\sigma z)}{\mathcal V_n(\sigma z)}>4\delta_n$,
and we still have
$\pair{\widetilde I_\lambda'(z)}{\mathcal V_n(z)}>4\delta_n$.

Choose a locally Lipschitz, $\sigma$-invariant
$\chi_n:\mathcal E\to[0,1]$ with $\chi_n(z)=1$ if
$\operatorname{dist}_{\mathcal E}(z,A_n)\le\frac32\delta_n$ and
$|\widetilde I_\lambda(z)-c|\le\varepsilon_n$, and $\chi_n(z)=0$ if
$\operatorname{dist}_{\mathcal E}(z,A_n)\ge2\delta_n$ or
$|\widetilde I_\lambda(z)-c|\ge2\varepsilon_n$. Since $\chi_n=0$ outside
$\mathcal N_n$, the field $-\chi_n\mathcal V_n$, extended by zero
outside $\mathcal O_n$, is locally Lipschitz and bounded on
$\mathcal E$. Let $\eta_n$ be its flow. Then
\begin{equation*}
 \begin{gathered}
 \eta_n(0,z)=z,\qquad
 \norm{\eta_n(t,z)-z}_{\mathcal E}\le t,\qquad
 \eta_n(t,\sigma z)=\sigma\eta_n(t,z),\\
 \frac{\mathrm d}{\mathrm dt}
 \widetilde I_\lambda(\eta_n(t,z))
 =-\chi_n(\eta_n(t,z))
    \pair{\widetilde I_\lambda'(\eta_n(t,z))}{\mathcal V_n(\eta_n(t,z))}
 \le0.
 \end{gathered}
\end{equation*}
Let $z\in A_n$ and $0\le t\le\delta_n$. Then
$\operatorname{dist}_{\mathcal E}(\eta_n(t,z),A_n)\le t<\frac32\delta_n$
and $\widetilde I_\lambda(\eta_n(t,z))\le c+\varepsilon_n$. So
$\chi_n(\eta_n(t,z))=1$ as long as
$\widetilde I_\lambda(\eta_n(t,z))\ge c-\varepsilon_n$. Then the energy
decreases at rate at least $4\delta_n$. Since
$4\delta_n^2=4\varepsilon_n$,
\begin{equation*}
 \widetilde I_\lambda(\eta_n(\delta_n,z))
 \le\max\bigl\{c-\varepsilon_n,\
       \widetilde I_\lambda(z)-4\delta_n^2\bigr\}
 \le c-\varepsilon_n \qquad(z\in A_n).
\end{equation*}
Since $\sup_{S^{k-1}}I_\lambda(\gamma)<0<c-2\varepsilon_n$, we have
$\chi_n=0$ at each $(0,\gamma(\omega))$. So the flow fixes these points.
Thus
\begin{equation*}
 \widehat h_n(y)=\eta_n(\delta_n,(0,h_n(y)))
       \in\widetilde\Gamma,\qquad
 \max_{y\in\overline B^k}
       \widetilde I_\lambda(\widehat h_n(y))
                      \le c-\varepsilon_n<c.
\end{equation*}
This contradicts \eqref{eq:equal-minimax} and proves
\eqref{eq:shadowing}.

\emph{Step 3: back to $X$.} All points of $A_n$ have first coordinate
zero. So
\begin{equation*}
 |\theta_n|\le\operatorname{dist}_{\mathcal E}(z_n,A_n)
          \le2\delta_n<1.
\end{equation*}
Let $u_n=T_{\theta_n}v_n$. By \eqref{eq:augmented-derivative} and
\eqref{eq:shadowing},
\begin{equation*}
 \begin{split}
 |I_\lambda(u_n)-c|&\le2\varepsilon_n,\qquad
 |P_\lambda(u_n)|\le8\delta_n,\\
 \norm{I_\lambda'(u_n)}_{X^*}
 &=\sup_{\norm{\varphi}_X\le1}
   \bigl|\pair{\partial_u\widetilde I_\lambda(\theta_n,v_n)}
                     {T_{-\theta_n}\varphi}\bigr|
 \le8\delta_n\norm{T_{-\theta_n}}_{\mathcal L(X)}
 \le C\delta_n\longrightarrow0.
 \end{split}
\end{equation*}
Thus $\{u_n\}$ satisfies \eqref{eq:PSP} at level $c$.
Lemma~\ref{lem:PSP} gives a critical point $u$ of $I_\lambda$ with
$I_\lambda(u)=c$.
\end{proof}

\begin{lemma}\label{lem:general-levels}
Assume \eqref{V2} and \eqref{g1}--\eqref{g4}. There exist numbers
$0<a_j\le b_j<\infty$ $(j\ge1)$, independent of
$\lambda\in(0,\lambda_0]$, with $a_j\to\infty$, such that for every
$j\ge1$ and every $\lambda\in(0,\lambda_0]$ there is
$u_{\lambda,j}\in X$ satisfying
\begin{equation*}
 I_\lambda'(u_{\lambda,j})=0,\qquad
       a_j\le I_\lambda(u_{\lambda,j})\le b_j.
\end{equation*}
\end{lemma}

\begin{proof}
Fix $q\in(2,2^*)$, and let $\tau_k$ be given by Lemma~\ref{lem:tau}. We
take $A=(M_0+1)/2$ in \eqref{eq:G-upper}. Since $V\ge-M_0$ and
$\lambda\ge0$, for every $\varepsilon>0$,
\begin{equation}\label{eq:general-lower}
 \begin{split}
 I_\lambda(u)
 &\ge\frac12|\nabla u|_2^2+\frac12\int_{\R^N}V_+(x)u^2\dx
       -\frac{M_0}2|u|_2^2+A|u|_2^2
       -C_\varepsilon|u|_q^q-\varepsilon|u|_{2^*}^{2^*}\\
 &=\frac12\norm{u}_V^2
       -C_\varepsilon|u|_q^q-\varepsilon|u|_{2^*}^{2^*}.
 \end{split}
\end{equation}
The constant $C_\varepsilon$ depends on $\varepsilon$. So we first fix
the radius and then choose the codimension. For $j\ge1$, let $R_j=j$,
and choose $\varepsilon_j>0$ with
$\varepsilon_j S^{-2^*/2}R_j^{2^*}\le R_j^2/8$. Then choose integers
$1\le k_1<k_2<\cdots$ with
$C_{\varepsilon_j}\tau_{k_j-1}^{q}R_j^q\le R_j^2/8$. Let
$u\in Z_{k_j-1}$ with $\norm{u}_V=R_j$. Then $|u|_q\le\tau_{k_j-1}R_j$,
and, by \eqref{eq:Sobolev},
$|u|_{2^*}^{2^*}\le S^{-2^*/2}|\nabla u|_2^{2^*}\le S^{-2^*/2}R_j^{2^*}$.
So \eqref{eq:general-lower} gives
\begin{equation}\label{eq:general-high-sphere}
 \inf_{\substack{u\in Z_{k_j-1}\\\norm{u}_V=R_j}}
         I_\lambda(u)\ge\frac12R_j^2-\frac18R_j^2-\frac18R_j^2
         =\frac14R_j^2=:a_j\longrightarrow+\infty,
\end{equation}
uniformly in $\lambda$.

Let $\gamma_j:S^{k_j-1}\to X$ be the odd map given by
Lemma~\ref{lem:negative-maps} with $k=k_j$ and $R_0=R_j$. Let
\begin{equation*}
 \begin{split}
 \Gamma_j&=\{h\in C(\overline B^{k_j},X):
       h\text{ is odd},\ h=\gamma_j\text{ on }S^{k_j-1}\},\\
 c_{\lambda,j}&=\inf_{h\in\Gamma_j}
             \max_{y\in\overline B^{k_j}}I_\lambda(h(y)).
 \end{split}
\end{equation*}
The map $h_j^0(t\omega)=t\gamma_j(\omega)$ ($0<t\le1$,
$\omega\in S^{k_j-1}$), $h_j^0(0)=0$, belongs to $\Gamma_j$. So
$\Gamma_j\ne\varnothing$. Let $h\in\Gamma_j$. Then
$\norm{h(\omega)}_V=\norm{\gamma_j(\omega)}_V>R_j$ on $S^{k_j-1}$. By
Lemma~\ref{lem:intersection} with $k=k_j$ and $\rho=R_j$, there is $y$
with $h(y)\in Z_{k_j-1}$ and $\norm{h(y)}_V=R_j$. Hence
$c_{\lambda,j}\ge a_j$ by \eqref{eq:general-high-sphere}. Also, since
$I_\lambda\le I_{\lambda_0}$,
\begin{equation*}
 c_{\lambda,j}\le
 \max_{y\in\overline B^{k_j}}I_{\lambda_0}(h_j^0(y))=:b_j<\infty.
\end{equation*}
By Lemma~\ref{lem:augmented-minimax}, there is a critical point
$u_{\lambda,j}$ of $I_\lambda$ with
$I_\lambda(u_{\lambda,j})=c_{\lambda,j}\in[a_j,b_j]$.
\end{proof}

\begin{proof}[Proof of Theorem~\ref{thm:general-multiplicity}]
Let $u_{\lambda,j}$, $a_j$ and $b_j$ be as in
Lemma~\ref{lem:general-levels}. Choose $j_1<j_2<\cdots$ with
$a_{j_{m+1}}>b_{j_m}$ for $m\ge1$. Fix $\lambda_n\to0^+$ with
$\lambda_n\le\lambda_0$. For each $m$, we apply
Theorem~\ref{thm:general-limit} to $\{u_{\lambda_n,j_m}\}_n$. Along a
subsequence, which may depend on $m$, we get a radial finite-energy weak
solution $u_m$ with
\begin{equation*}
 u_{\lambda_n,j_m}\to u_m\text{ in }H^1(\R^N),\qquad
 a_{j_m}\le I_0(u_m)=\lim_{n\to\infty}
       I_{\lambda_n}(u_{\lambda_n,j_m})\le b_{j_m}.
\end{equation*}
By Lemma~\ref{lem:Pohozaev} with $\lambda=0$, $g(u_m)u_m\in L^1(\R^N)$.
Since $a_{j_{m+1}}>b_{j_m}$ and $a_{j_m}\to\infty$,
\begin{equation*}
 0<I_0(u_1)<I_0(u_2)<\cdots,\qquad I_0(u_m)\longrightarrow\infty.
\end{equation*}
\end{proof}

\ack{C. Huang was supported by the National Natural Science Foundation
of China (Grant No. 12601203) and YNNUMA (No. 2510). }

\data{No new data were created or analysed in this study.}

\end{document}